   \documentclass[final,3p,sort&compress]{elsarticle}

\usepackage{amssymb,amsmath,amsthm,amsbsy}
\usepackage{physics}
\usepackage[pdftex,bookmarksnumbered]{hyperref}
\usepackage{graphicx}
\usepackage{epstopdf}
\usepackage{epsfig}
\usepackage{multicol}
\usepackage{subfigure}
\usepackage{multirow}
\usepackage{algorithm,algorithmic}
\usepackage{bm}
\usepackage{setspace}
\usepackage{color}
\usepackage{lineno}

\journal{Elsevier}

\begin{document}


\begin{frontmatter}



\title{Multi-fidelity regression using artificial neural networks: efficient approximation of parameter-dependent output quantities}



\author[add1]{Mengwu Guo\corref{cor1}}
\ead{m.guo@utwente.nl}
\cortext[cor1]{Corresponding author.}

\author[add2]{Andrea Manzoni}
\ead{andrea1.manzoni@polimi.it}

\author[add3]{Maurice Amendt}

\author[add2]{Paolo Conti}
\ead{paolo5.conti@mail.polimi.it}

\author[add3]{Jan S. Hesthaven}
\ead{Jan.Hesthaven@epfl.ch}

\address[add1]{Department of Applied Mathematics, University of Twente}
\address[add2]{MOX -- Dipartimento di Matematica, Politecnico di Milano}
\address[add3]{Institute of Mathematics, \'{E}cole Polytechnique F\'{e}d\'{e}rale de Lausanne}

\begin{abstract}
   Highly accurate numerical or physical experiments are often very time-consuming or expensive to obtain. When time or budget restrictions prohibit the generation of additional data, the amount of available samples may be too limited to provide satisfactory model results. Multi-fidelity methods deal with such problems by incorporating information from other sources, which are ideally well-correlated with the high-fidelity data, but can be obtained at a lower cost. By leveraging correlations between different data sets, multi-fidelity methods often yield superior generalization when compared to models based solely on a small amount of high-fidelity data. In the current work, we present the use of artificial neural networks applied to multi-fidelity regression problems. By elaborating a few existing approaches, we propose new neural network architectures for multi-fidelity regression. The introduced models are compared against a traditional multi-fidelity regression scheme -- co-kriging. A collection of artificial benchmarks are presented to measure the performance of the analyzed models. The results show that cross-validation in combination with Bayesian optimization consistently leads to neural network models that outperform the co-kriging scheme. Additionally, we show an application of multi-fidelity regression to an engineering problem. The propagation of a pressure wave into an acoustic horn with parametrized shape and frequency is considered, and the index of reflection intensity is approximated using the proposed multi-fidelity models. A finite element, full-order model and a reduced-order model built through the reduced basis method are adopted as the high- and low-fidelity, respectively. It is shown that the multi-fidelity neural network returns outputs that achieve a comparable accuracy to those from the expensive, full-order model, using only very few full-order evaluations combined with a larger amount of inaccurate but cheap evaluations of a reduced order model.

\end{abstract}

\end{frontmatter}


\section{Introduction}

Artificial neural networks (ANNs) have arguably been one of the most active topics during the recent years. They have been successfully applied in a substantial number of research areas, including image recognition \cite{ImageRec}, translation \cite{Chen_2019}, and fraud detection \cite{Ghosh1994CreditCF}. More recently, ANNs have also been widely used in the emerging area of machine learning in computational science and engineering, sometimes referred to as scientific machine learning \cite{baker2019workshop}. The remarkable expressive power of neural networks (NNs) has made them stand out in the solution of forward and inverse problems governed by partial differential equations (PDEs) \cite{raissi2019physics,weinan2017deep,sirignano2018dgm}, reduced order modeling \cite{hesthaven2018non,lee2020model,fresca2020comprehensive}, data-driven discovery \cite{ray2018artificial}, multiscale analysis \cite{regazzoni2020machine}, and so on. This success can largely be explained by three major factors: computational power, flexibility and access to large data sets. The great flexibility of NNs, as well as their multi-purpose nature, is a dominant factor explaining their overwhelming success, not only in research, but also in real-world applications. Several application-dependent architectures have been developed, e.g., convolutional neural networks are often used in image related tasks, whereas recurrent neural networks have found their success in speech recognition. Nevertheless, even for a fixed structure, a neural network is still able to adapt to different situations. This feature is mainly explained by the large number of parameters and hyperparameters that can be tuned to fit data in many situations.

The ability to handle large data sets without overwhelming computational costs has made NNs a good candidate for multi-fidelity (MF) regression. MF regression exploits correlations between different data sets to provide a regression model that generalizes better than a simple regression model, which only takes a single data set into account. The most common setup for MF regression is a set of data sources with different fidelity levels. High-fidelity (HF) samples are usually rare due to their cost, be it computational or experimental. However, such data are accurate and are the best available knowledge about the problem. In contrast, low-fidelity (LF) data are assumed to be easy and cheap to obtain, yet might lack accuracy. Often generated numerically, the LF data set ideally expresses major trends of the problem and correlate well with the HF data. The goal of MF regression is to infer trends from LF data and use it to approximate the HF model, especially in regions where HF data are sparse.

A natural MF setting arose in the geostatistics community, who realized that \textit{"if core data at other locations are correlated, they should be included to improve the regression"} \cite{geostat}. In the context of soil porosity, precise measurements are combined with seismic data to better predict porosity in large areas. Another common MF situation arises when solving PDEs \cite{raissi2017inferring,kast2020non}. High order numerical schemes with a fine mesh give rise to accurate, HF data, whereas the usage of a coarse mesh, a partially converged solution, or a linearized equation can lead to LF data. In the past decade, MF methods have found applications in many areas of scientific computing, including uncertainty quantification, inference, and optimization. We refer to \cite{peherstorfer2018survey} for a comprehensive review. A widely used MF technique is co-kriging \cite{Hagan2000,alvarez2012kernels} which relies on vector-valued Gaussian processes for regression. Such a Gaussian process regression scheme presents two major benefits. Firstly, as a non-parametric regression tool, it is suitable for many different applications. Secondly, the method can be cast in the Bayesian framework, and the regression results naturally include an uncertainty estimation, which is usually desirable. Nevertheless, Gaussian process regression is not suitable for many applications as it suffers from different drawbacks, such as the curse of dimensionality. Consequently, new methods for MF regression are needed, ideally able to detect non-trivial, highly nonlinear correlations between the data sets of different fidelity levels, and should be applicable to a general class of problems. 

Based on these considerations, ANNs appear to be a promising candidate to solve MF problems. In the related current work, we consider MF regression with ANNs. Several approaches have been proposed in the literature so far and successfully used NNs in a MF setting. In \cite{Meng}, the authors tested a deep NN structure on different artificial MF benchmarks, extended the idea to the physics-informed NN scheme and applied it to inverse problems governed by partial differential equations (PDEs). Their MF model clearly outperforms a single fidelity regression. A different NN architecture, incorporated into a Monte Carlo sampling algorithm to estimate uncertainties in a MF setting, was introduced in \cite{Motamed}, and reduced computational cost has been observed as compared to traditional Monte Carlo sampling. In addition, different MF strategies for training NNs were discussed in \cite{aydin2019general}, NNs were used to approximate the discrepancy between the HF and LF physics-constrained NNs in \cite{liu2019multi}, and deep neural networks (DNNs) were embedded into co-kriging in \cite{raissi2016deep}. Moreover, a composition of Gaussian processes in a multi-layer network structure was used for MF modeling in \cite{cutajar2019deep}, and a multi-fidelity Bayesian neural network scheme has been developed in \cite{meng2020multi} and applied to the physics-informed versions.

In this work, we propose different ANN architectures for the purpose of MF regression. Inspired by \cite{Meng} and \cite{Motamed}, we present two all-in-one models in which different fidelity levels are trained simultaneously, as well as two multilevel models that define separate NNs for the hierarchy of fidelity levels. In addition, we utilize a strategy based on cross-validation and Bayesian optimization to automatically select the best performing NN hyperparameters. To the best of our knowledge, hyperparameter optimization (HPO) has not been investigated for MF neural networks. The goal of this work is to define a reliable strategy that consistently proposes a neural network which achieves high accuracy and shows good generalization properties. We will assess the performance of the proposed NNs on a set of manufactured benchmarks, all chosen carefully to test for the desirable properties. All models will be compared to single-fidelity regression schemes. Comparisons will also include co-kriging results to benchmark the proposed models against common practices in MF frameworks. 

Additionally, we show an application of the proposed MF regression schemes to the evaluation of a quantity of interest defined as a functional of the solution to a parameter-dependent problem governed by PDEs. Such a task often occurs in the applied sciences and engineering, where multiple evaluations of PDE solutions, for different scenarios described in terms of physical or geometrical parameters, can be computationally demanding if relying on high-fidelity, full-order models (FOMs) such as detailed finite element approximations. To overcome this difficulty, low-fidelity, reduced-order models (ROMs) can be built through the reduced basis (RB) method. Despite ROMs featuring much lower-dimensional solution spaces than those of the FOMs, they are able to capture the critical physical features of the FOMs. A reduced model seeks the solutions on a low-dimensional manifold which is approximated by a linear trial subspace spanned by a set of global basis functions, built from a set of full-order snapshots. The ROM accuracy is often granted at the price of a relatively large number of basis functions involved in the reduced-order approximation. On the other hand, an efficient assembly of the ROM during the online stage may only be possible provided that an expensive hyper-reduction is performed during the offline stage. Therefore, a low dimensionality without expensive hyper-reduction can make a reduced model extremely efficient, but potentially inaccurate. Our goal, enabled by the MF neural network schemes in this work, is to provide accurate approximations to the output quantities by leveraging a relatively large number of output evaluations using very low-dimensional ROMs and a small number of evaluations using the FOM, so as to avoid the efficiency issues stemming from the ROM construction and evaluation without compromising the outcome accuracy. In particular, we apply the proposed approaches to a parametrized PDE problem, namely the propagation of a pressure wave into an acoustic horn with parametrized shape and frequency, described by the Helmholtz equation. We assess the impact of both the quality and the amount of the training data on the overall accuracy of the MF regression outcome.

Following the introduction, the concepts of ANNs and Gaussian process regression are briefly reviewed, and their underlying correlation is discussed in Section 2. Several NN structures for MF regression are introduced and discussed in Section 3, and their effectiveness is demonstrated by a series of benchmark test cases in Section 4. An application to a parametrized PDE problem is presented in Section 5, and conclusions are drawn in Section 6.

\section{Artificial neural networks and Gaussian processes for regression}\label{ANN}

\subsection{Artificial neural networks (ANNs)}

In this section, we consider an $L$-hidden-layer fully-connected NN \cite{strang2019linear} with hidden layers of width $N_l$ for the $l$-th layer and the nonlinear activation function $\phi$, $1\leq l \leq L$. At the $j$-th neuron in the $l$-th layer of the NN, the pre- and post-activation are denoted by $z_j^l$ and $x_j^l$, respectively, $1\leq i \leq M_l$, $M_l$ being the width of the $l$-th layer. Let $\vb{x}=\vb{x}^0\in\mathbb{R}^{d_\text{in}}$ denote the inputs of the network and $\vb{y}=\vb{z}^{L+1}\in\mathbb{R}^{d_\text{out}}$ denote the outputs. Note that we have $M_0=d_\text{in}$ and $M_{L+1}=d_\text{out}$. Weight and bias parameters between the $(l-1)$-th and $l$-th layers are represented by $W_{ij}^{l}$ and $b_i^l$, respectively, $1\leq l\leq (L+1)$, $1\leq i \leq M_{l}$, $1\leq j \leq M_{l-1}$. Then one has
\begin{equation}
\begin{split}
 & z_i^l(\vb{x})=b_i^l + \sum_{j=1}^{M_{l-1}}W_{ij}^l x_j^{l-1}(\vb{x})\,,\quad
 x_i^{l}(\vb{x})=\phi(z_i^{l}(\vb{x}))\,,\quad 1\leq i \leq M_l, ~1\leq l \leq L\,, \quad \text{and}\\
 & y_i(\vb{x})=z_i^{L+1}(\vb{x})=b_i^{L+1} + \sum_{j=1}^{M_{L}}W_{ij}^{L+1} x_j^{L}(\vb{x})\,,\quad \ \ \, \qquad  1\leq i \leq d_\text{out}\,.
 \end{split}
\end{equation}

A multivariate function $\vb{y}=\vb{f}(\vb{x})$ is approximated by a vector-valued network surrogate $\vb{f}^\texttt{NN}(\cdot;\vb{W},\vb{b}):\mathbb{R}^{d_\text{in}}\to \mathbb{R}^{d_\text{out}}$ to be trained on the input-output pairs $\{(\vb{x}^{(k)},\vb{y}^{(k)})\}_{k=1}^{N}$. Here $\vb{W}$ and $\vb{b}$ are the vectors collecting all the weight and bias parameters, respectively, and $N$ is the number of data pairs. Such a training is often performed by minimizing a cost function:
\begin{equation}
(\vb{W},\vb{b})=\arg\min_{\vb{W},\vb{b}} \left\{ \frac{1}{N}\sum_{k=1}^{N}\|\vb{y}^{(k)}-\vb{f}^\texttt{NN}(\vb{x}^{(k)};\vb{W},\vb{b})\|_2^2 + \lambda \|\vb{W}\|_2^2\right\}\,,
\end{equation}
in which the first term is the mean square error (MSE) and the second term is a regularization term with $\lambda \geq 0$ being the penalty coefficient.

\subsection{Gaussian process regression (GPR)}

\subsubsection*{Single-fidelity GPR}

A Gaussian process (GP) is a collection of random variables, any finite number of which obeys a joint Gaussian distribution. In the GPR, the prior on the scalar-valued regression function $f:\mathbb{R}^{d_\text{in}}\to \mathbb{R}$ is assumed to be a GP corrupted by an independent Gaussian noise term, i.e., for $(\vb{x},\vb{x}')\in \mathbb{R}^{d_\text{in}} \times \mathbb{R}^{d_\text{in}}$,
\begin{equation}
f(\vb{x})\sim \text{GP}(0,\kappa(\vb{x},\vb{x}'))\,,\quad y=f(\vb{x})+\epsilon\,,\quad \epsilon\sim\mathcal{N}(0,\chi^2)\,,
\end{equation}
where $\chi$ is the standard deviation of a Gaussian noise term $\epsilon$, and the semi-positive definite kernel function $\kappa$ gives the covariance of the prior GP. 

Given $N$ pairs of input-output training data, a prior joint Gaussian is defined for the corresponding outputs:
\begin{equation}
\vb*{y}|\vb*{X}~\sim~\mathcal{N}(\vb{0},\vb{K}_y)\,,\quad \vb{K}_y=\mathrm{Cov}[\vb*{y}|\vb*{X}]=\kappa(\vb*{X},\vb*{X})+\chi^2\vb{I}_N\,,
\end{equation}
where $\vb*{y}=\{y^{(1)},y^{(2)},\cdots,y^{(N)}\}^{\mathrm{T}}$, $\vb*{X}=[~\vb{x}^{(1)}~|~\vb{x}^{(2)}~|~\cdots~|~\vb{x}^{(N)}~]$ and $\vb{I}_N$ is the $M$-dimensional unit matrix, $N$ being the number of training samples.

From a regression model, the goal is to predict the noise-free output $f^*(\vb{s})$ for a new test input $\vb{s}\in\mathbb{R}^{d_\text{in}}$. By the standard rules for conditioning Gaussians, the posterior predictive distribution conditioning on the training data is obtained as a new GP:
\begin{equation}
\begin{split}
& f^*(\vb{s})|~\vb{s},\vb*{X},\vb*{y}~\sim~\text{GP}(m^*(\vb{s}),c^*(\vb{s},\vb{s}'))\,, \\
& m^*(\vb{s}) = \kappa(\vb{s},\vb*{X})\vb{K}_y^{-1}\vb*{y}\,,\quad c^*(\vb{s},\vb{s}') = \kappa(\vb{s},\vb{s}')-\kappa(\vb{s},\vb*{X})\vb{K}_y^{-1}\kappa(\vb*{X},\vb{s}')\,.
\end{split}
\end{equation}

\subsubsection*{Multi-fidelity GPR}

GPR with training data from different fidelity levels is known as cokriging \cite{Hagan2000} or vector-valued GPR \cite{alvarez2012kernels}.  In such a regression scheme, one can use a large amount of LF data and only a limited number of HF samples to training a model of a reasonable accuracy. Since the LF evaluations are cheap, the cost of training data preparation can be reduced by controlling the number of HF evaluations. Assuming a linear correlation between the different fidelity levels, we can employ the linear model of coregionalization (LMC) \cite{alvarez2012kernels} that expresses the prior of a hierarchy of $D$ solution fidelities as
\begin{equation}\label{vvgpr}
f_i (\vb{x})= \sum_{j=1}^D a_{i,j} u_j(\vb{x}) \,, \quad i = 1,2,\cdots, D\,,
\end{equation}
i.e., each level of solution $f_i$ is written as a linear combination of $D$ independent Gaussian processes $u_j \sim \text{GP}(0,\kappa_j(\cdot,\cdot))$. In addition, the vector $\vb{a}_j$, $1\leq j \leq D$, collects the weights of the corresponding GP component $u_j$, i.e., $\vb{a}_j=\{a_{1,j},\cdots,a_{D,j}\}^\text{T}$. This formulation leads to a matrix-valued kernel for the MF GPR model as
\begin{equation}
\vb*{\mathcal{K}}(\vb{x},\vb{x}') = \sum_{j=1}^{D}\vb{a}_j\vb{a}_j^\text{T} \kappa_j(\vb{x},\vb{x}')\,.
\end{equation}

In the two-level case, the well-known form of AR(1)-cokriging \cite{Hagan2000} is a special form of the linear model  defining the prior of a LF solution $f_L$ and a HF $f_H$:
\begin{equation}
\begin{split}
& f_H(\vb{x})= \rho u_1(\vb{x})+ u_2(\vb{x})\,,\\
& f_L(\vb{x}) = u_1(\vb{x})\,,
\end{split}
\end{equation}
in which $\vb{a}_1=\{\rho,1\}^\text{T}$ and $\vb{a}_2=\{1,0\}^\text{T}$. Conditioning on the training input-output pairs from both the LF and HF, denoted by $(\vb*{X}_L,\vb*{y}_L)$ and $(\vb*{X}_H,\vb*{y}_H)$, respectively, the predictive distribution for the HF can be expressed as a posterior GP, i.e., $f_H^*(\vb{s}) |~\vb{s},\vb*{X}_H,\vb*{y}_H,\vb*{X}_L,\vb*{y}_L \sim \text{GP}$.

\subsection{The link between ANNs and GPR}

It can be shown that the prior of a neural network output can be seen as a set of Gaussian processes under the following probabilistic assumptions \cite{lee2017deep,neal2012bayesian}\footnote{The discussion here involves some modification from the work in \cite{lee2017deep,neal2012bayesian}.}: (I) All the weight parameters $W_{ij}^l$'s are independent and identically distributed (i.i.d), as are all the bias parameters $b_i^{l}$'s, and the weight and bias parameter sets are independent of each other; (II) In the $l$-th layer, $1\leq 2 \leq L+1$, $b_i^l \sim \mathcal{N}(0,\sigma_b^2)$, and $W_{ij}^l$'s are independently drawn from any distribution with zero mean and variance $\sigma_w^2 / M_{l-1}$; and (III) $M_l \to \infty$, $1\leq l \leq L$. Here we show by induction that $\{z^{l}_i:1\leq i \leq M_{l}\}$ are i.i.d. zero-mean Gaussian processes and  $\{x^{l}_j:1\leq j \leq M_{l}\}$ are i.i.d. for all $2\leq l \leq L+1$. 

We  consider an arbitrary set of finite locations of the input $\vb{x}$, denoted by $\vb*{X}$. Since both $\{b_i^1:\forall i\}$ and $\{W_{i,j}^1:\forall i,j\}$ are i.i.d., one obtains $\{z_i^1(\vb*{X})=b_i^1+\sum_{j=1}^{d_\text{in}}W_{ij}^1 x_j(\vb*{X}):\forall i\}$. Thus $\{x_j^1(\vb*{X})=\phi(z_j^1(\vb*{X})): \forall j\}$ are i.i.d. Then it can be recovered that $\{z_i^2(\vb*{X})=b_i^2+\sum_{j=1}^{M_1}W_{ij}^2 x_j^1(\vb*{X}):\forall i\}$ are i.i.d. For each $1\leq i \leq M_2$, we apply the multivariate central limit theorem \cite{durrett2019probability} to a sequence of i.i.d. random vectors $\{\sqrt{M_1}W^2_{i1}x_1^1(\vb*{X}),\sqrt{M_1}W^2_{i2}x_2^1(\vb*{X}),\cdots\}$, all with zero mean and covariance $\sigma_w^2 \text{Cov}[x_\cdot^1(\vb*{X})]$, and we obtain that $z_i^2(\vb*{X})\sim \mathcal{N}(\vb{0},\sigma_b^2\vb{I}+\sigma_w^2 \mathbb{E}[x_\cdot^1(\vb*{X})\otimes x_\cdot^1(\vb*{X})])$ as $M_1 \to \infty$. Since the input locations $\vb*{X}$ are arbitrary, each $z_i^2(\vb{x})$ is a Gaussian process as $z_i^2(\cdot)\sim \text{GP}(0,K^2(\cdot,\cdot))$ with its kernel $K^1$ defined as $K^2(\vb{x},\vb{x}')=\sigma_b^2\vb{I}+\sigma_w^2 \mathbb{E}[\phi(z_\cdot^1(\vb{x}))\phi(z_\cdot^1(\vb{x}'))]$. After the nonlinear activation in the 2nd layer, we have the i.i.d. $\{x^{2}_j:\forall j\}$. Therefore the proposition holds true for $l=2$.

Assume it holds true for $l$, $2\leq l \leq L$. The proposition can be verified to be true for $l+1$, similarly to that for $l=2$. The pre-activation in each intermediate layer follows a Gaussian process $z_\cdot^l\sim\text{GP}(0,K^l(\cdot,\cdot))$, $2\leq l \leq L+1$, and 
\begin{equation}\label{rckernel}
K^{l}(\vb{x},\vb{x}')= \sigma_b^2 +\sigma_w^2 \mathbb{E}_{z_{\cdot}^{l-1}\sim \text{GP}(0,K^{l-1})}[\phi(z_\cdot^{l-1}(\vb{x}))\phi(z_\cdot^{l-1}(\vb{x}'))]\,.
\end{equation}
The outputs $\vb{y}(\vb{x})=\vb{z}^{L+1}(\vb{x})$ thus follow i.i.d. Gaussian process priors and the kernel function can be formed from the recursive relation \eqref{rckernel}.

\section{Artificial neural networks for multi-fidelity regression}\label{NNMFR}

Currently, there have been two successful approaches \cite{Meng,Motamed} to the use of NNs in a multi-fidelity context. The neural networks for multi-fidelity regression (NNMFR) presented in these two studies
show inherent differences at the architectural level of the networks. Whereas \cite{Meng} uses a single NN to perform a simultaneous regression of both the HF and LF data, \cite{Motamed} splits them up and used two distinct ANNs, one for each fidelity level. Note that in the present work, we restrict ourselves to bi-fidelity problems, i.e., we are interested in approximating the scalar HF function $f_\texttt{HF}(\vb{x})$ by incorporating information from the scalar LF function $f_\texttt{LF}(\vb{x})$. Consequently, there will generally be only a few available observations of the HF function, whereas the LF data are abundant. In that respect, it is useful to introduce the following notation:
\begin{itemize}
	\item $N_\texttt{HF}$ and $N_\texttt{LF}$ denote the number of HF and LF samples, respectively.
	\item The training set for the HF function is given by $\mathcal{T}_\texttt{HF} = \{(\vb{x}_\texttt{HF}^{(i)},y_\texttt{HF}^{(i)}):1\leq i \leq N_\texttt{HF}\}$, which satisfies the HF function $y_\texttt{HF}=f_\texttt{HF}(\vb{x})$.
	\item By replacing the subscripts $\texttt{HF}$ in the above notations by $\texttt{LF}$, we obtain the LF counterparts.
	\item The NN approximations corresponding to the HF and LF functions are denoted by $f_\texttt{HF}^\texttt{NN}(\vb{x})$ and $f_\texttt{LF}^\texttt{NN}(\vb{x})$, respectively. 
\end{itemize}
In the current section, we will review existing NNMFR approaches and propose our own strategies by elaborating and adapting existing ideas. Different from the existing NNMFR strategies, our all-in-one models consider the LF outputs as latent variables of the HF surrogate or mimic the correlation between HF and LF levels in co-kriging, and our multilevel models are formulated directly from the hierarchy of fidelity levels.

\subsection{All-in-one NNMFR} \label{oneforall}
As previously noted, a single NN, which performs two simultaneous regressions on the HF and LF data, is used in \cite{Meng}. It essentially resembles the structure of an autoencoder, where the inputs are encoded towards the LF output $y_\texttt{LF}$, then decoded and encoded again for the HF output $y_\texttt{HF}$. The setup relies on the assumption that
\begin{equation}\label{decompose}
f_\texttt{HF}(\vb{x}) = \mathcal{F}(\vb{x},y_\texttt{LF}) = \alpha\mathcal{F}_{l}(\vb{x},y_\texttt{LF})+ (1-\alpha)\mathcal{F}_{nl}(\vb{x},y_\texttt{LF}), \quad \alpha \in [0,1], 
\end{equation}
i.e., the unknown mapping $\mathcal{F}$ of the LF to the HF data can be decomposed into a linear and a nonlinear part, denoted by $\mathcal{F}_{l}$ and $\mathcal{F}_{nl}$, respectively.

In \eqref{decompose}, the hyperparameter $\alpha$ determines the strength of the linear correlation, with $\alpha = 1$ corresponding to a fully linear relation between the HF and LF outputs. The determination of the value of
 $\alpha$ was not specified in \cite{Meng}, neither are any values indicated in the reported results.

\begin{figure}[h]
	\centering
	\subfigure[Intermediate]
	{\includegraphics[width=0.48\linewidth]{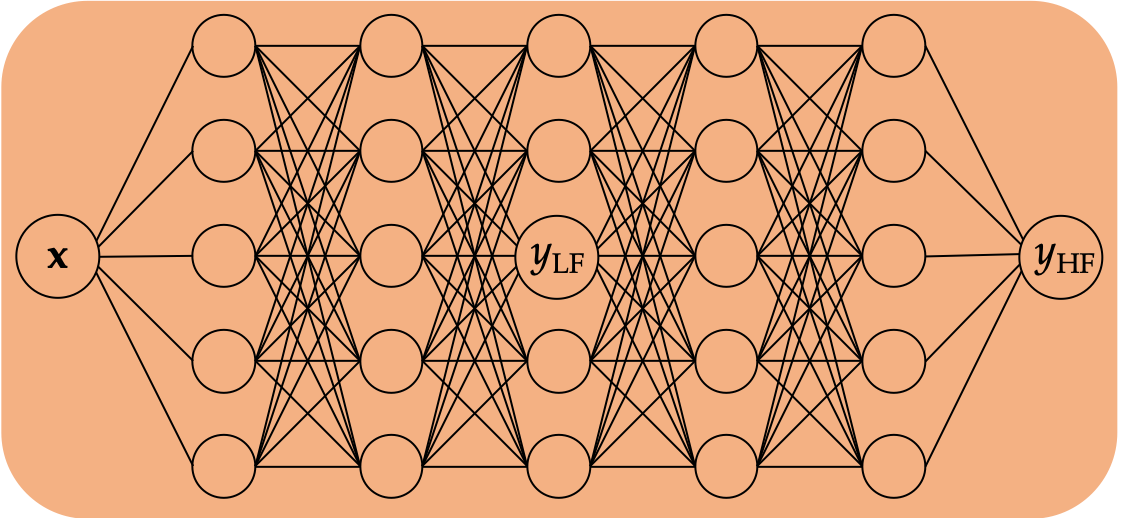}}
	\quad
	\subfigure[GPmimic]
	{\includegraphics[width=0.48\linewidth]{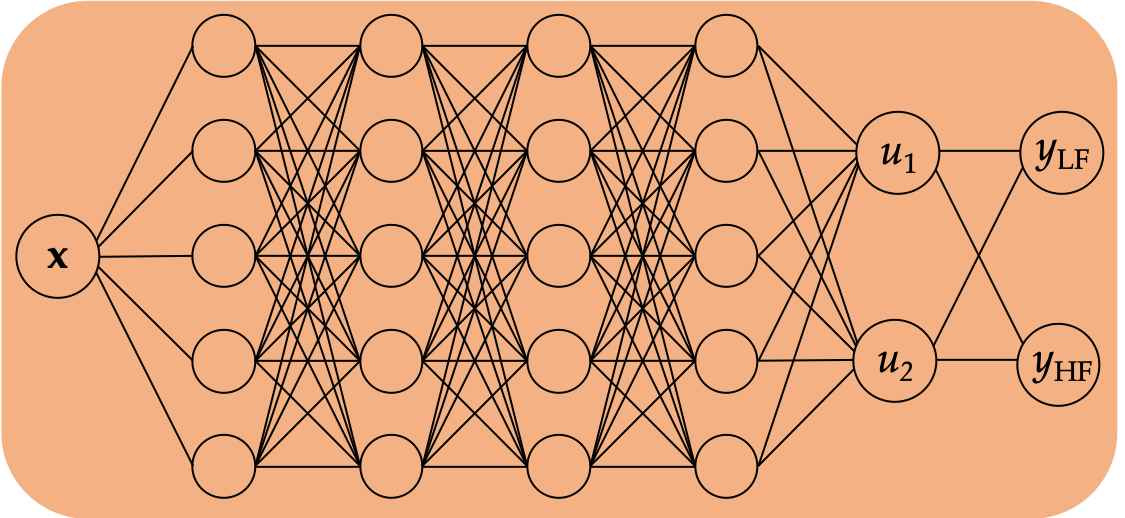}}
	\subfigure[2-step]
	{\includegraphics[width=0.40\linewidth]{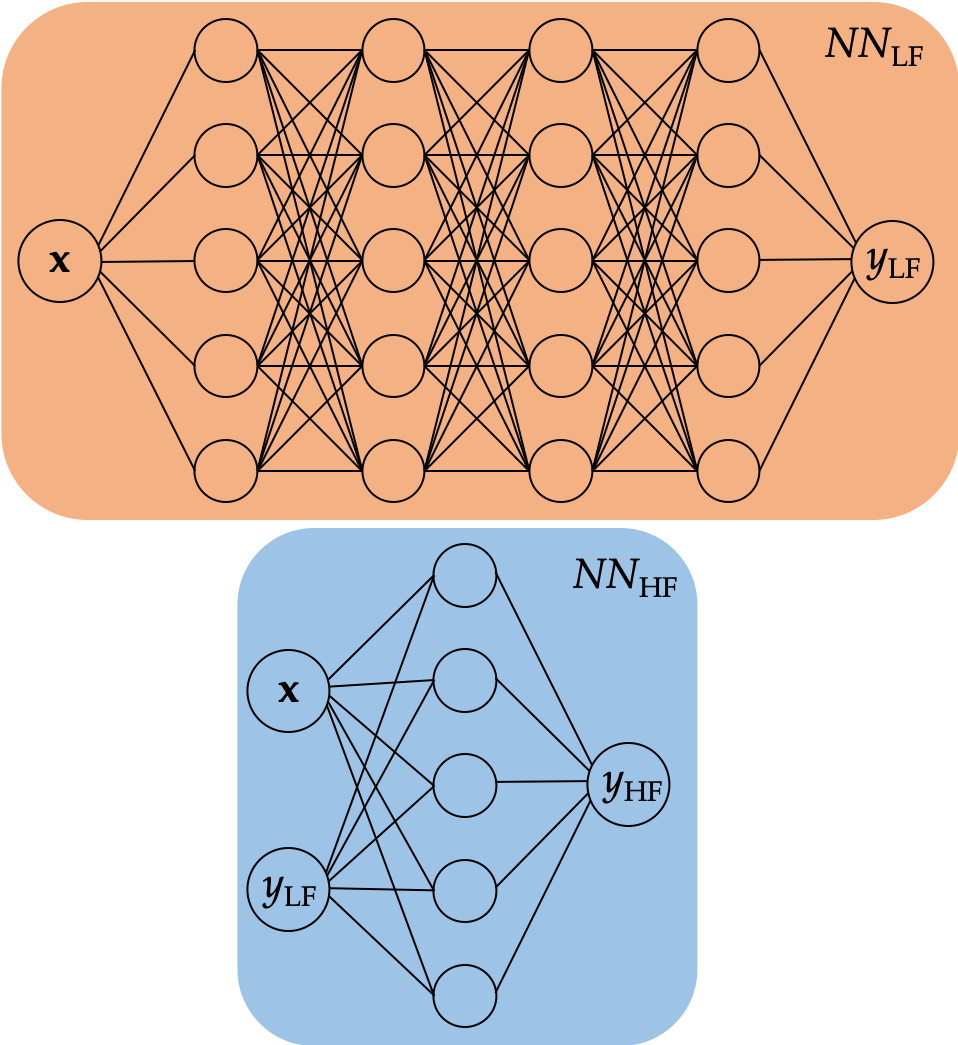}}
	\quad\quad
	\subfigure[3-step]
	{\includegraphics[width=0.40\linewidth]{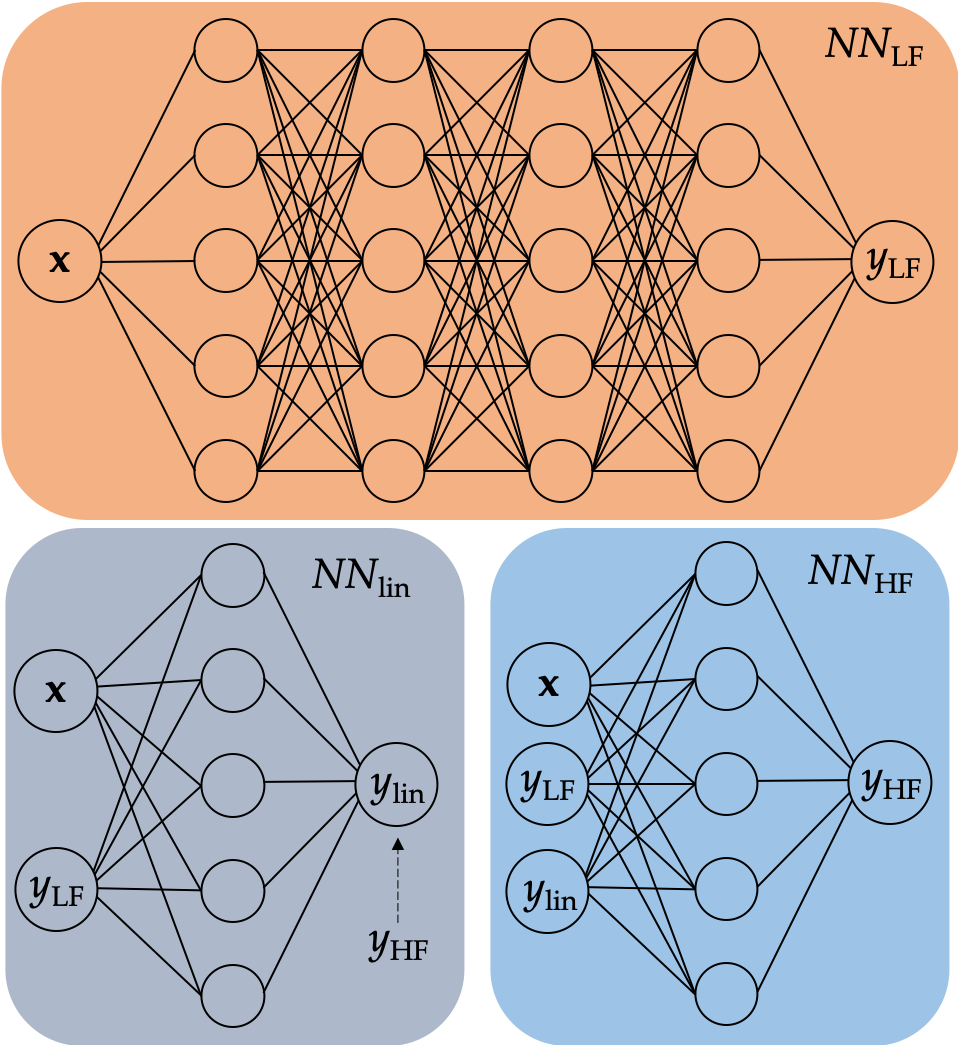}}
	\caption{Different ANNs proposed for MF regression. All-in-one models in the top row (a, b), and multilevel models in the bottom row (c, d). The 'Intermediate model' in (a) considers the LF outputs at intermediate latent variables of the surrogate for HF. The 'GPmimic model' in (b) employs the LMC to mimic MF GPR. The '2-step' and '3-step' models in (c) and (d) define seperate NNs for the hierarchy of fidelity levels. }
	\label{fig:nnmfr}
\end{figure}

Based on these considerations, we propose the following two modified NN architectures, see (a) and (b) in Fig. \ref{fig:nnmfr}:
\begin{enumerate}
	\item \textbf{"Intermediate" model}: The first NN structure is similar to the one in \cite{Meng}, except that the same input layer is used for HF and LF data and we omit the autoencoder resemblance by adding additional nodes to the layer containing the LF output. Furthermore, we do not impose \eqref{decompose} as we seek to directly model the function $\mathcal{F}(\vb{x},y_\texttt{LF})$. We refer to this NN as the "Intermediate" model, referring to the location of the LF output. In the numerical experiments, the Intermediate network has 5 hidden layers; the LF output is situated in the 3rd hidden layer and the HF output in the last layer. The number of neurons in the first 3 layers is fixed to 64 each and the widths of the layers between the LF and HF outputs are determined by the hyperparameter optimization. Except for the output, we use the hyperbolic tangent activation function.
	\item \textbf{"GPmimic" model}: The second NNMFR is based on the similarity between a wide NN and a GP. Based on this correspondence, and since Gaussian processes are commonly used for MF regression, we implement an architecture, which seeks to mimic the action of a GP. Hereafter we refer to this architecture as "GPmimic". In contrast to NNs, the use of GPs becomes infeasible in higher dimensions. Hence, the GPmimic NN indirectly enables GP-like regression even in the case of large data sets and high-dimensional inputs.
	The NN architecture is depicted in Fig. \ref{fig:nnmfr}. The neurons labelled by $u_1$ and $u_2$ are defined to be analogous to two independent Gaussian processes, especially when the previous layers are wide,  and they play the same roles as in the vector-valued GPR \eqref{vvgpr}. By considering no nonlinear activation in the output layer, the NN outputs of the GPmimic model are given as:
	\begin{equation}\label{GPmimicMF} 
	\begin{split}
	y_\texttt{HF}(\vb{x}) &= W_{11}u_1(\vb{x}) + W_{12}u_2(\vb{x})+ b_1\\
	y_\texttt{LF}(\vb{x}) &= W_{21}u_1(\vb{x}) + W_{22}u_2(\vb{x})+ b_2. 
	\end{split}
	\end{equation}
	By inspecting \eqref{GPmimicMF}, we notice that the outputs are recovered as an affine transformation of the two variables $u_1(\vb{x})$ and $u_2(\vb{x})$, and the parameters $W_{ij}'s$ define the correlations between $y_\texttt{HF}$ and $y_\texttt{LF}$ as in the vector-valued GPR. In principle, this setting only allows to model linear correlations between the HF and LF data.
	
\end{enumerate}

Remember that the all-in-one architectures perform vector-valued learning of $\vb{f} = [f_\texttt{HF}(\vb{x}),f_\texttt{LF}(\vb{x})]^\text{T}$. There are thus two different error components, one for each fidelity level. More specifically, these components are given by the training errors MSE$_\texttt{HF} $ and MSE$_\texttt{LF} $ of the HF and LF models, defined as
\begin{equation}
\text{MSE}_\texttt{HF} = \frac{1}{N_\texttt{HF}}\sum_{i=1}^{N_\texttt{HF}}|y_\texttt{HF}^{(i)} - f_\texttt{HF}^\texttt{NN}(\vb{x}_\texttt{HF}^{(i)})|^2, \quad \text{and} \quad
\text{MSE}_\texttt{LF}  = \frac{1}{N_\texttt{LF}}\sum_{i=1}^{N_\texttt{LF}}|y_\texttt{LF}^{(i)} - f_\texttt{LF}^\texttt{NN}(\vb{x}_\texttt{LF}^{(i)})|^2.
\end{equation}
Hence, training an all-in-one network involves a multiobjective minimization problem, which is generally not easy to solve. In \cite{Meng}, the problem is reduced to a single optimization problem by considering a loss function given by the addition of the MSEs related to the two fidelity sets and a $L_2$-regularization term.

A more general approach in multiobjective optimization is to use a weighted sum instead, to account for objectives of different scales. Even though the data sets can be initially scaled to ensure the same order of magnitude, large discrepancies between MSE$_\texttt{HF}$ and MSE$_\texttt{LF}$ could develop during the learning process. Based on these considerations, the loss function in the Intermediate and GPmimic models is given as
\begin{equation}
\mathcal{L} = \alpha \text{MSE}_\texttt{HF} + (1-\alpha) \text{MSE}_\texttt{LF} + \lambda \|\vb{W}\|_2^2,
\end{equation}
where $\alpha \in [0,1]$ (unrelated to the one in \eqref{decompose}) acts as a scaling factor between the two fidelity levels and $\lambda > 0$ is a penalty parameter. The extreme cases $\alpha = 0$ and $\alpha = 1$ correspond to the single-fidelity regressions for the LF and HF data, respectively. The choice $\alpha = 0.5$ represents a balanced importance of the LF and HF levels. The value of $\alpha$ will be tuned by hyperparameter optimization in the numerical examples. The $\alpha$ values minimize the testing errors of $k$-fold cross-validation and are determined by Bayesian optimization.

\subsection{Multilevel NNMFR}

The NNMFR introduced in Subsection \ref{oneforall} relies on a single NN to learn the multidimensional function $\mathbf{f} = [f_\texttt{HF}(\vb{x}),f_\texttt{LF}(\vb{x})]^\text{T}$. Another approach, leading to the multilevel NNMFR, is to use distinct NNs to model these functions. Such an approach has been employed  in \cite{Motamed}, and the presented method can be summarized as follows. A first NN $NN_1$ learns a correlation function $y_\texttt{HF}=\mathcal{F}(\vb{x},y_\texttt{LF})$ based on the input data $\left\{(\vb{x}^{(i)}_\texttt{HF},y_\texttt{LF}(\vb{x}^{(i)}_\texttt{HF})):1\leq i \leq N_\texttt{HF}\right\}$ and the output data $\left\{y^{(i)}_\texttt{HF}: 1\leq i \leq N_\texttt{HF}\right\}$. In this case, it is required that the HF input locations should form a subset of the LF input locations, i.e., $\mathcal{X}_\texttt{HF}=\{\vb{x}_\texttt{HF}^{(i)}: 1\leq i \leq N_\texttt{HF}\} \subseteq \mathcal{X}_\texttt{LF}=\{\vb{x}_\texttt{LF}^{(i)}: 1\leq i \leq N_\texttt{LF}\} $. For each available LF sample $\vb{x}^{\prime} \in \mathcal{X}_\texttt{LF} \setminus \mathcal{X}_\texttt{HF}$, an approximate HF data sample $y'_\texttt{HF}(\vb{x}')$ can be generated by the network $NN_1$, and these generated data, together with the existing HF data at $\mathcal{X}_\texttt{HF}$, can be used as the new HF data set $\mathcal{Y}'_\texttt{HF}=\{y_\texttt{HF}(\mathcal{X}_\texttt{HF}),y'_\texttt{HF}(\mathcal{X}_\texttt{LF} \setminus \mathcal{X}_\texttt{HF})\}$ . Finally, a second NN $NN_2$ is trained to model the HF function $f_\texttt{HF}(\vb{x})$ based the input-output pairs between $\mathcal{X}_\texttt{LF}$ and $\mathcal{Y}'_\texttt{HF}$.
Following this multilevel approach in \cite{Motamed}, we propose two modified multi-step NNs for MF regression by adopting the following major changes:
\begin{enumerate}
	\item The LF function $f_\texttt{LF}(\vb{x})$ is modeled by a first NN $NN_\texttt{LF}$. This modification is relevant if the computational cost of additional LF data is not cheap, for instance when having to solve PDEs.
	\item Now that we can rapidly generate new LF data using $NN_\texttt{LF}$, modeling the direct mapping between $\vb{x}$ and $y_\texttt{HF}=f_\texttt{HF}(\vb{x})$ is unnecessary. Hence our second network $NN_\texttt{HF}$ will be analogous to $NN_1$ in \cite{Motamed}.
\end{enumerate}
These considerations lead to the following multilevel architectures, also see (c) and (d) in Fig. \ref{fig:nnmfr}:
\begin{enumerate}
	\item \textbf{"2-step" model}: A DNN $NN_\texttt{LF}$ is trained on $\mathcal{T}_\texttt{LF}$ to learn the LF function $f_\texttt{LF}(\vb{x})$. Using the NN $NN_\texttt{LF}$, we can predict the values of the LF function at the training inputs $\mathcal{X}_\texttt{HF}$ of the HF data, denoted by $f_\texttt{LF}^\texttt{NN}(\mathcal{X}_\texttt{HF})=\{f_\texttt{LF}^\texttt{NN}(\vb{x}_\texttt{HF}^{(i)}): 1 \leq i \leq N_\texttt{HF}\}$. Then a second artificial network $NN_\texttt{HF}$ approximates the HF function $y_\texttt{HF}=\mathcal{F}(\vb{x},y_\texttt{LF})$ based on the input data $(\mathcal{X}_\texttt{HF},f_\texttt{LF}^\texttt{NN}(\mathcal{X}_\texttt{HF}))=\{(\vb{x}_\texttt{HF}^{(i)},f_\texttt{LF}^\texttt{NN}(\vb{x}_\texttt{HF}^{(i)})): 1 \leq i \leq N_\texttt{HF}\}$ and the available HF output data $\mathcal{Y}_\texttt{HF}=y_\texttt{HF}(\mathcal{X}_\texttt{HF})=\{y_\texttt{HF}^{(i)}: 1 \leq i \leq N_\texttt{HF}\}$. $NN_\texttt{HF}$ is a shallow NN consisting of a single hidden layer.
	
	In other words, the network $NN_\texttt{HF}$ approximates the function $y_\texttt{HF}=\mathcal{F}(\vb{x},y_\texttt{LF})$ based on the HF and LF data at the same locations $\mathcal{X}_\texttt{HF}$. However, as the HF and LF data locations, i.e., $\mathcal{X}_\texttt{HF}$ and $\mathcal{X}_\texttt{LF}$, are generated independently, the observations of the LF function at the HF inputs $\mathcal{X}_\texttt{HF}$ have to be evaluated from $NN_\texttt{LF}$ if not directly available at $\mathcal{X}_\texttt{HF}$.

	\item \textbf{"3-step" model}: This model is a modification of the 2-step model by adding an additional level of fidelity, generated by a third network $NN_{\text{lin}}$. $NN_\text{lin}$ is equivalent to the $NN_\texttt{HF}$ in the 2-step model in the sense that it models the correlation function $y_\texttt{HF}=\mathcal{F}(\vb{x},y_\texttt{LF})$. However, no nonlinear activation function is used in $NN_{\text{lin}}$. In other words, $NN_{\text{lin}}$ is responsible for capturing the linear correlations between the HF and LF data sets. Let $f_\text{lin}^\texttt{NN}\left(\mathcal{X}_\texttt{HF},f_\texttt{LF}^\texttt{NN}(\mathcal{X}_\texttt{HF})\right)$ denote the set of outputs of $NN_\text{lin}$ at the HF input locations $\mathcal{X}_\texttt{HF}$. This set will then serve as part of the inputs for the third and final network $NN_\texttt{HF}$. $NN_\texttt{HF}$ approximates the correlation function $y_\texttt{HF}=\mathcal{F}'(\vb{x},y_\texttt{LF},y_\text{lin})$, $y_\text{lin}$ being the output of $NN_\text{lin}$ as a linear approximation of $y_\texttt{HF}$, and the training of $NN_\texttt{HF}$ is based on the input-output pairs between $\left(\mathcal{X}_\texttt{HF},f_\texttt{LF}^\texttt{NN}(\mathcal{X}_\texttt{HF}),f_\text{lin}^\texttt{NN}\left(\mathcal{X}_\texttt{HF},f_\texttt{LF}^\texttt{NN}(\mathcal{X}_\texttt{HF})\right)\right)$ and $\mathcal{Y}_\texttt{HF}$.
	$NN_\texttt{HF}$ is again a shallow NN consisting of a single hidden layer. Note that the 3-step model should only present an advantage over the 2-step model when there exist strong linear correlations between the HF and LF data.
\end{enumerate}

\vspace{3mm}
It should be clear that the proposed architectures present a few differences. The key difference between the all-in-one and the multilevel strategies is at the structural level. All-in-one NNMFR aims to approximate the two-dimensional function $\mathbf{f} = [f_\texttt{HF}(\vb{x}),f_\texttt{LF}(\vb{x})]^\text{T}$ with a single NN, while multilevel NNMFR uses distinct networks to model the HF and LF functions separately. However, to incorporate information from the LF function, this is done sequentially by first modeling $f_\texttt{LF}(\vb{x})$ and then the correlation between the two fidelity levels.

We note that all-in-one models include the additional hyperparameter $\alpha$. This hyperparameter can potentially be used to incorporate prior knowledge into the regression model. Suppose for instance that we know that the LF data provide a good representation of the problem. We can then enforce the NN to accurately model the LF data by setting the value of $\alpha$ accordingly. Furthermore, this can be taken into account in the Intermediate model by choosing the number of neurons in the layer containing the LF output, where a smaller layer width may indicate a stronger dependency on the LF data. In the multilevel architectures, the question of how much trust is put in the LF data is left to the model and determined during the training process. 


An advantage of the multilevel approaches is that the LF function does not have to be approximated by a neural network. Instead, we can use other techniques which might be more accurate or less time-consuming. For instance, in the presence of discontinuities, an accurate approximation can be obtained by utilizing gradient boosting algorithms \cite{gradientboost}. 

Finally, in the non-hierarchical cases where the available information sources present similar levels of fidelities, the GPmimic model is the only NNMFR that can be used without any adaptation. In fact, all other proposed structures rely on a hierarchy of the available data sets according to their fidelity levels.  

\section{Numerical results (I): benchmark test cases}

In this section, we analyze the performance of the proposed ANN structures in MF regression problems. For that purpose, the models presented in Section \ref{NNMFR} are tested on a variety of artificial benchmarks, and   compared against commonly used co-kriging methods.

In all benchmarks, the hyperparameters of the neural networks, for instance the parameter $\alpha$ to balance the MSE terms, are optimized by $k$-fold cross-validation and Bayesian optimization, whereas the hyperparameters of the GPs are obtained by maximizing the marginal likelihood \cite{Rasmussen2004}. For the Bayesian optimization, Python package Hyperopt \cite{hyperopt} is used for a tree-structured parzen estimator (TPE) based approach \cite{bergstra2011algorithms}, in which the priors on hyperparameters can be chosen from a range of distributions for both continuous and discrete random variables. Both the single-fidelity and multi-fidelity GPRs are implemented by the Python package GPy \cite{gpy2014}. Furthermore, GPy automatically adds a noise term for each fidelity level to account for noisy data. As artificial benchmarks give rise to noiseless observations, the standard deviation values of all the noise terms are fixed to $10^{-5}$ in these benchmarks.

The numerical study in this section compares the performance of the proposed NNs against co-kriging. In each case, the qualitative difference between single- and multi-fidelity regressions is noted briefly before presenting a comparative study between NNMFR and multi-fidelity GPR. Models are evaluated using the MSE on a test set covering the whole input domain. The results will reveal whether choosing a model based on the validation error is a good strategy. The comparison between NN regression and the best GPR will not only include accuracy results but also computational time. In addition, we evaluate the following $R^2$ score that allows for a cross-benchmark comparison:
\begin{equation}
R^2 = 1 - \dfrac{\sum_{i=1}^{N_\texttt{test}}\left(y_\texttt{HF}^{(i)} - f^{(i)}_\texttt{reg}\right)^2}{\sum_{i=1}^{N_\texttt{test}}\left(y_\texttt{HF}^{(i)} - \bar{y}_\texttt{HF}\right)^2}\,,\quad  \text{ where } \bar{y}_\texttt{HF} = \dfrac{1}{N_\texttt{test}}\sum_{i=1}^{N_\texttt{test}}y^{(i)}_\texttt{HF}\,,
\end{equation}
$N_\texttt{text}$ is the size of the test set, and $f_\texttt{reg}$ denotes the prediction given by a regression model.


\subsection{Benchmark case 1: Linear correlation} \label{1dcase}

\begin{figure}[h!]
	\centering
    \subfigure[SFR with GP, MFR with co-kriging.]{\includegraphics[width=0.475\textwidth]{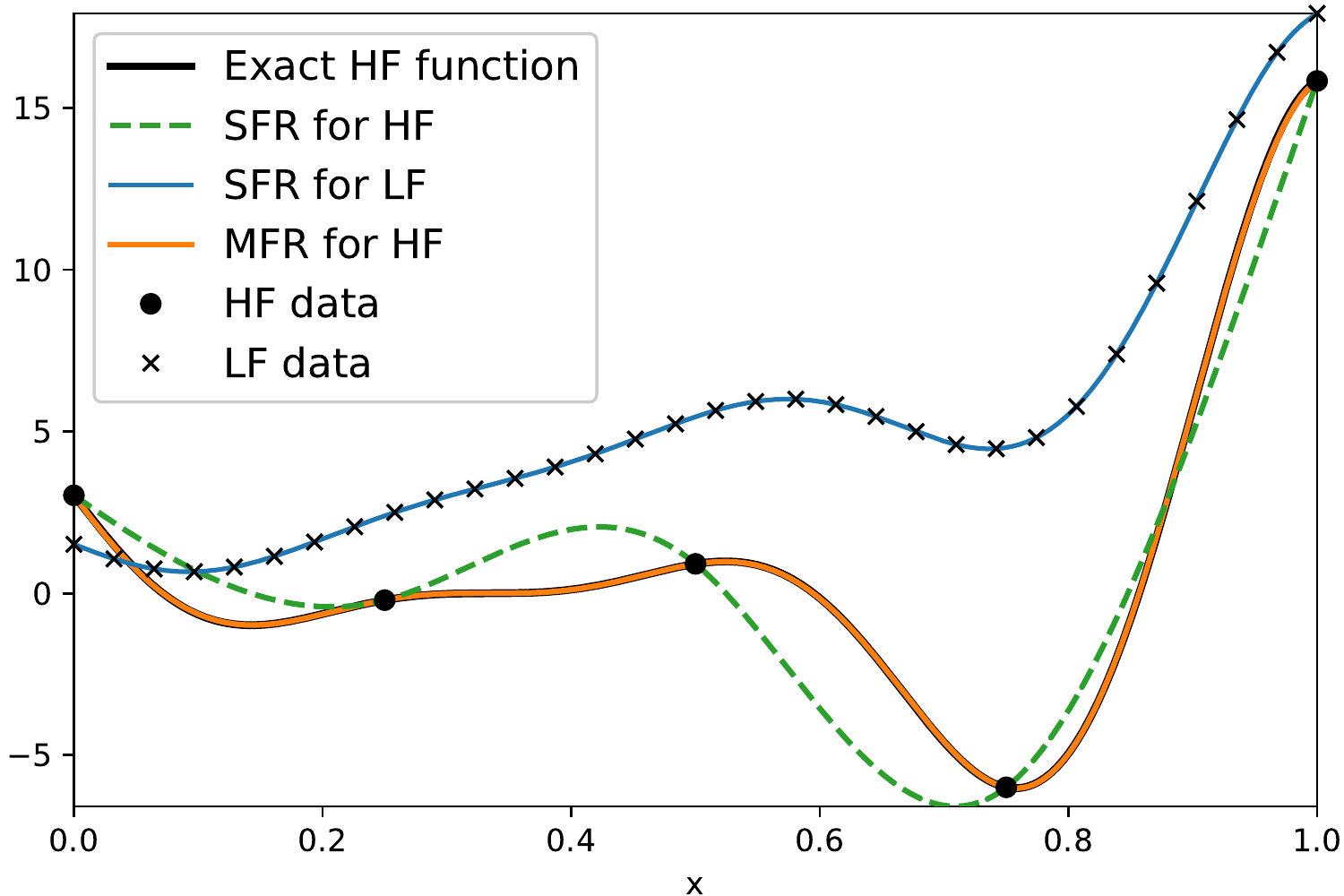}}
	\quad
	\subfigure[SFR with NN, MFR with 3-step.]{\includegraphics[width=0.475\textwidth]{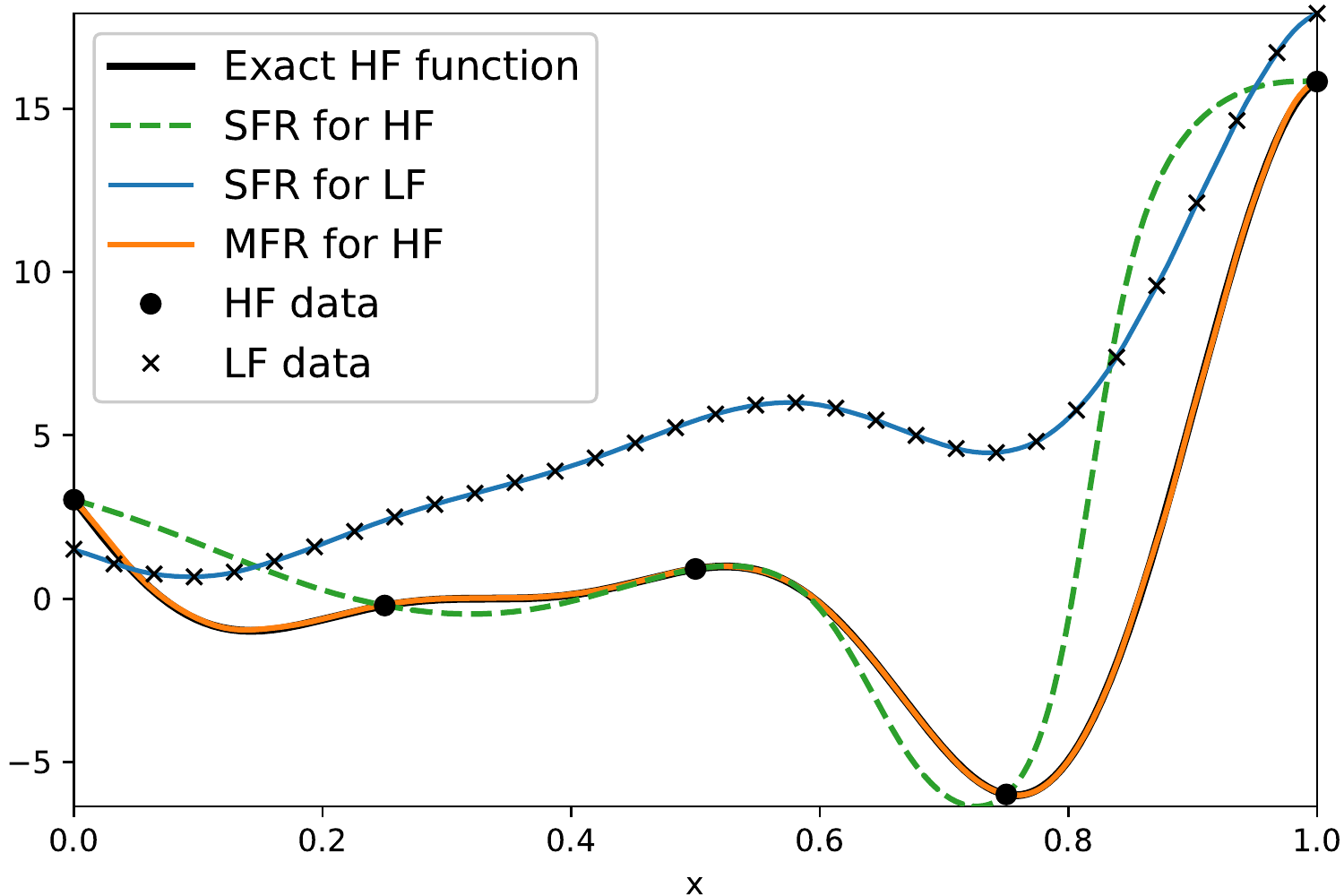}}
	\par\bigskip
	\subfigure[SFR with NN, MFR with 2-step.]{\includegraphics[width=0.475\textwidth]{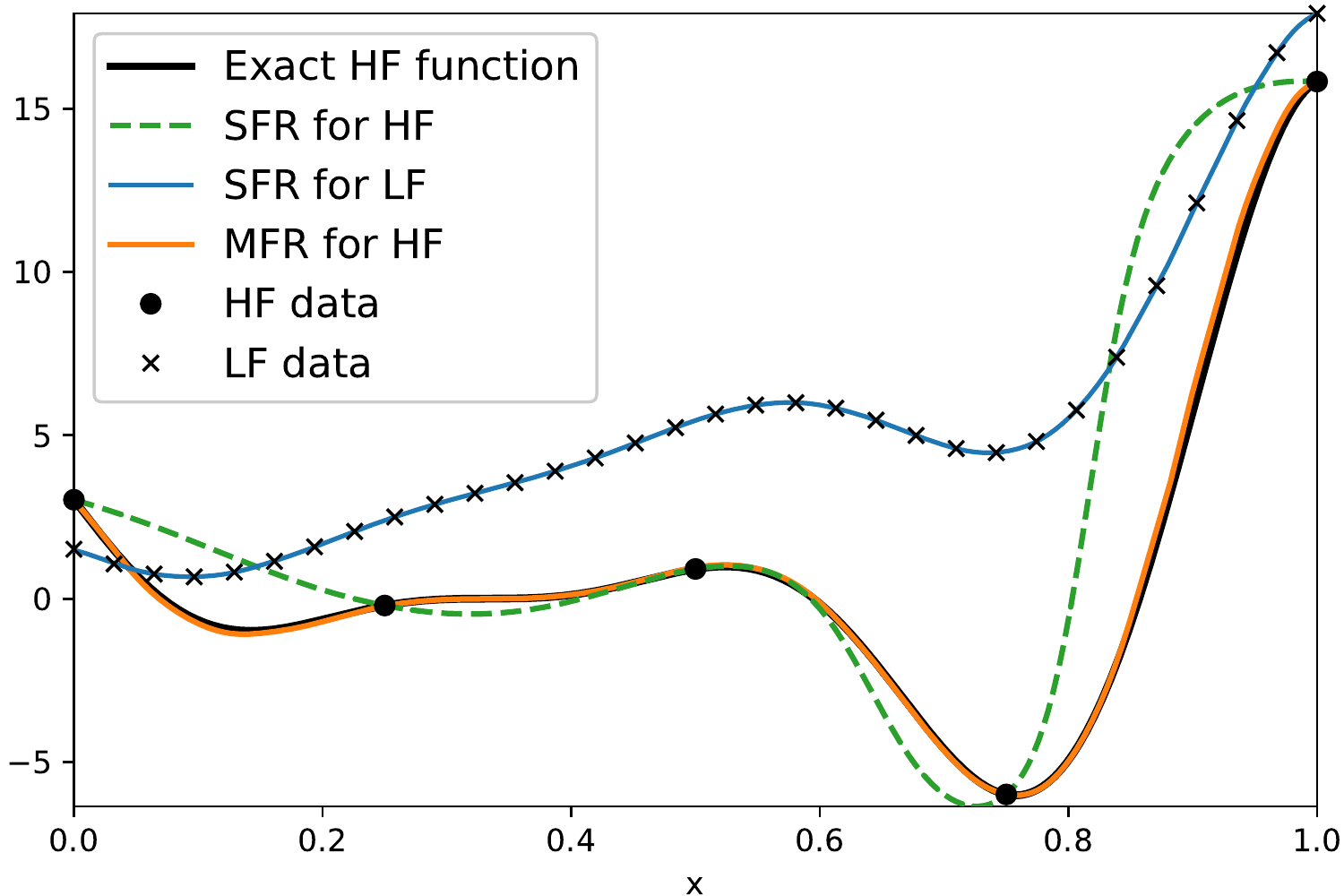}}
	\quad
	\subfigure[SFR with NN, MFR with GPmimic.]{\includegraphics[width=0.475\textwidth]{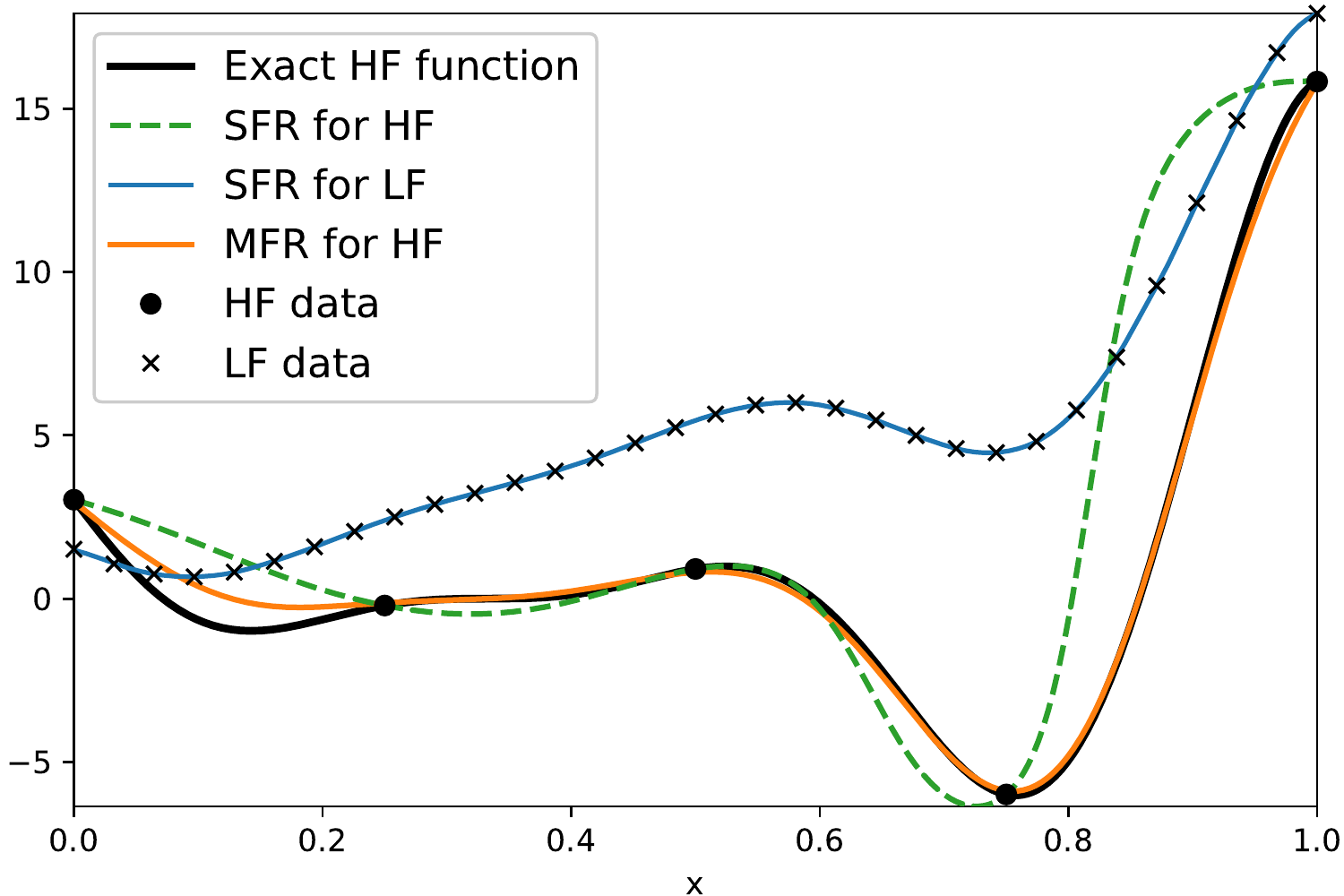}}
	\par\bigskip
	\subfigure[SFR with NN, MFR with Intermediate.]{\includegraphics[width=0.475\textwidth]{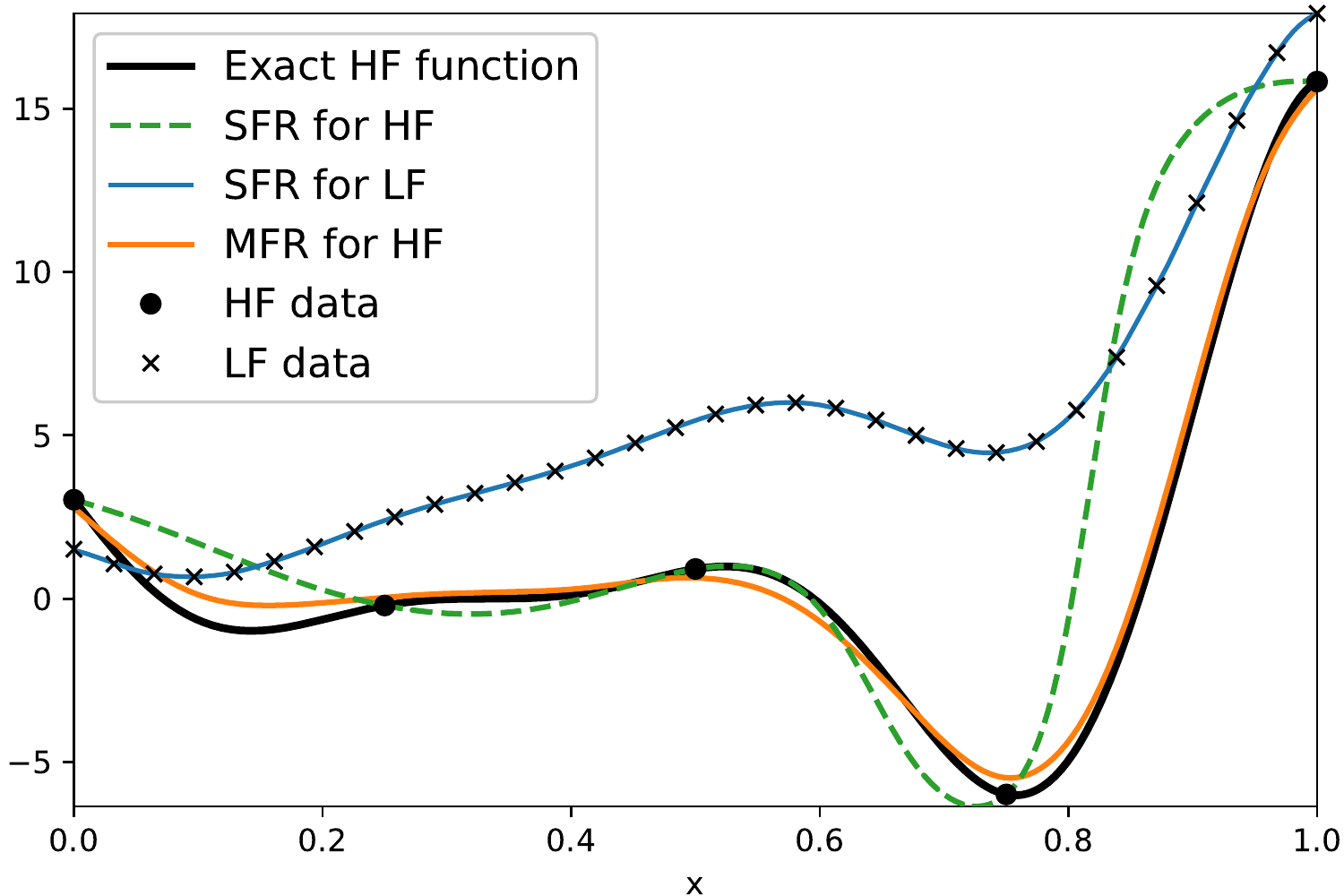}}
	\caption{Linear correlation: single- and MF regression results using GPs and different ANNs. There are 5 HF (circle) and 32 LF (cross) observations. Single-fidelity model based on the 5 HF data points is shown in green (dashed).}
	\label{fig:bench1results}
\end{figure}

The first benchmark is a common test case for MF methods. HF and LF functions are defined over $\Omega = [0,1]$ as
\begin{align*}
f_\texttt{HF}(x) &=  (6x-2)^2 \sin(12x-4), \\
f_\texttt{LF}(x) &= 0.5f_\texttt{HF}(x) + 10(x-0.5) + 5,
\end{align*}
respectively. The LF function is obtained by a linear transformation from the HF function and the input variable $x$. The setup serves as a good initial test for MF models. HF and LF samples are given by 5 and 32 equally spaced values in $\Omega$, respectively. Since the number of the HF data is very limited, hyperparameters of the NNs are optimized using leave-one-out cross-validation (LOOCV). The values can be found in Table \ref{summary}. Fig. \ref{fig:bench1results} shows the results of both the single-fidelity regression (SFR) and MF regression (MFR). Both the NN and the GPR fail to approximate the function based solely on the HF data; in particular, both models fail to accurately predict the function in the interval $[0,0.8]$. It is also important to note that the kernel used for the GPR is induced by an NN structure \cite{Meng,lee2017deep}. In fact, all other kernels predict an almost constant mean function with large uncertainty intervals. Table \ref{bench1summary} provides quantitative insight in the regression results. No validation error is available for co-kriging, as the hyperparameters are optimized by maximizing the marginal likelihood. However, since the best performing kernel has been chosen based on the MSE on the test set, it is fair to assume that the performance of GPR is slightly overestimated in general.

In general, except for the 2-step model, the validation errors are very conservative estimates of the test error. Co-kriging and the 3-step model outperform all other models by one or even two orders of magnitude. This can be explained by their mathematical setup that is designed to exactly model the linear correlation between the HF and LF functions. In fact, the 3rd neural network $NN_\texttt{HF}$ of the 3-step model is redundant in this case, as there is no nonlinear trend to catch. Consequently, in this first and simple benchmark there is nothing to be gained by using ANN, as they all perform worse than classic GPR and are also computationally more expensive.

\begin{table}[h!]
	\caption{Linear correlation: comparison of the MF regression models. Indicated times account for both HPO and final model predictions.}
	\renewcommand{\arraystretch}{1.3}
	\centering%
	\begin{tabular}{c|cccc}
		\hline%
		Model & Validation MSE & Test MSE & $R^2$ & Elapsed time (s)\\
		\hline
		Co-kriging & -  & $\mathbf{2.9\times10^{-4}}$  & 0.999& 9\\
		Intermediate & 3.03$\times10^{1}$  & $1.87\times10^{-1}$  & 0.990 & 650\\
		GP-mimic & 4.25$\times10^{0}$ & 1.55$\times10^{-1}$ &0.992 & 752\\
		2-step & 7.81$\times10^{-3}$ & 3.37$\times10^{-2}$ & 0.998 & 104\\
		3-step & 2.23$\times10^{-3}$ & 9.53$\times10^{-4}$ & 0.999 & 401\\
		\hline
	\end{tabular}
	\label{bench1summary}
\end{table}

\subsection{Benchmark case 2: Discontinuous function}

\begin{figure}[h!]
	\centering
	\subfigure[HF and LF functions.]{\includegraphics[width=.475\textwidth]{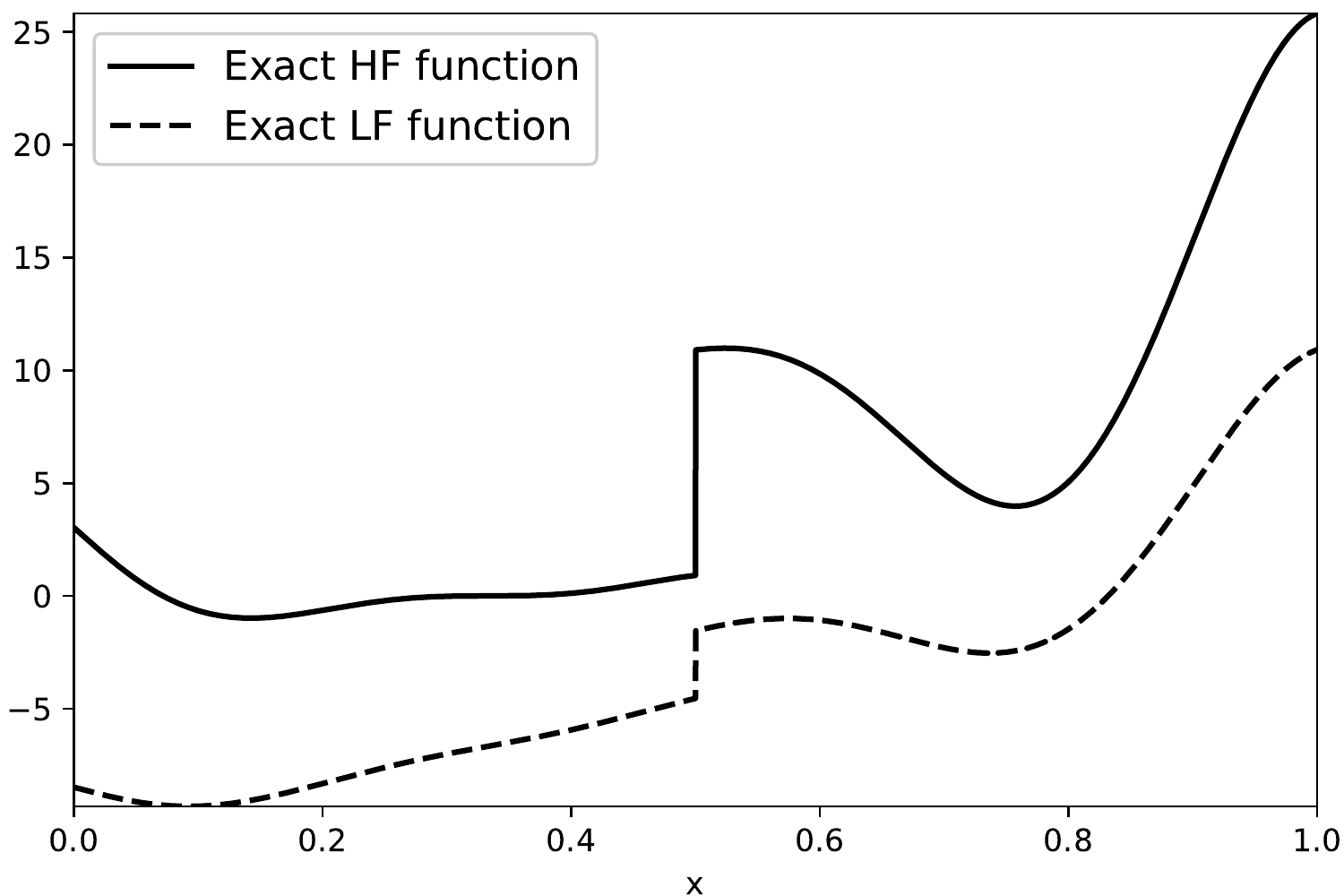}\label{fig:bench2_functions}}
	\quad
	\subfigure[SFR with GP, MFR with co-kriging.]{\includegraphics[width=.475\textwidth]{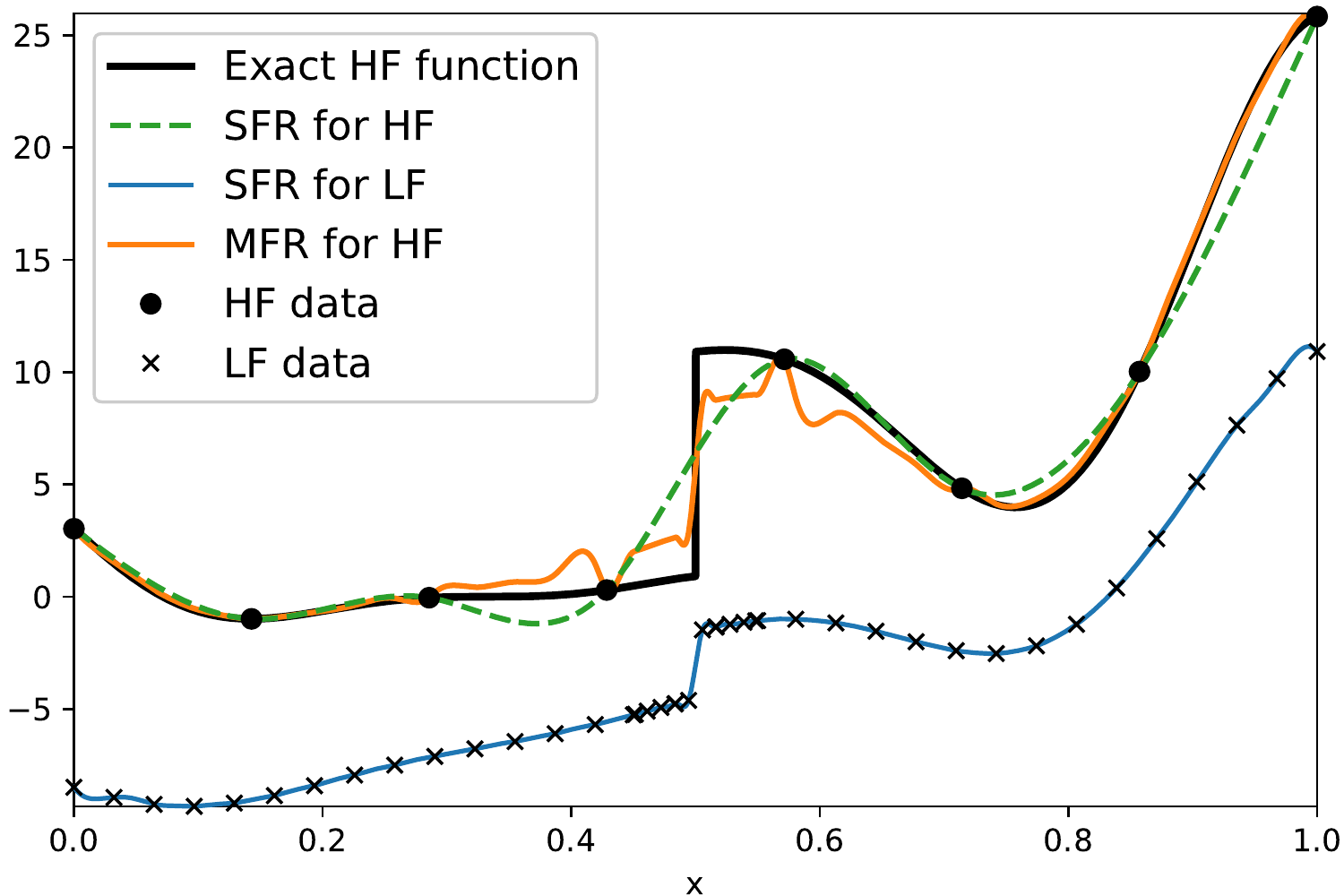}}
	\par\bigskip
	\subfigure[SFR with NN, MFR with 2-step.]{\includegraphics[width=.475\textwidth]{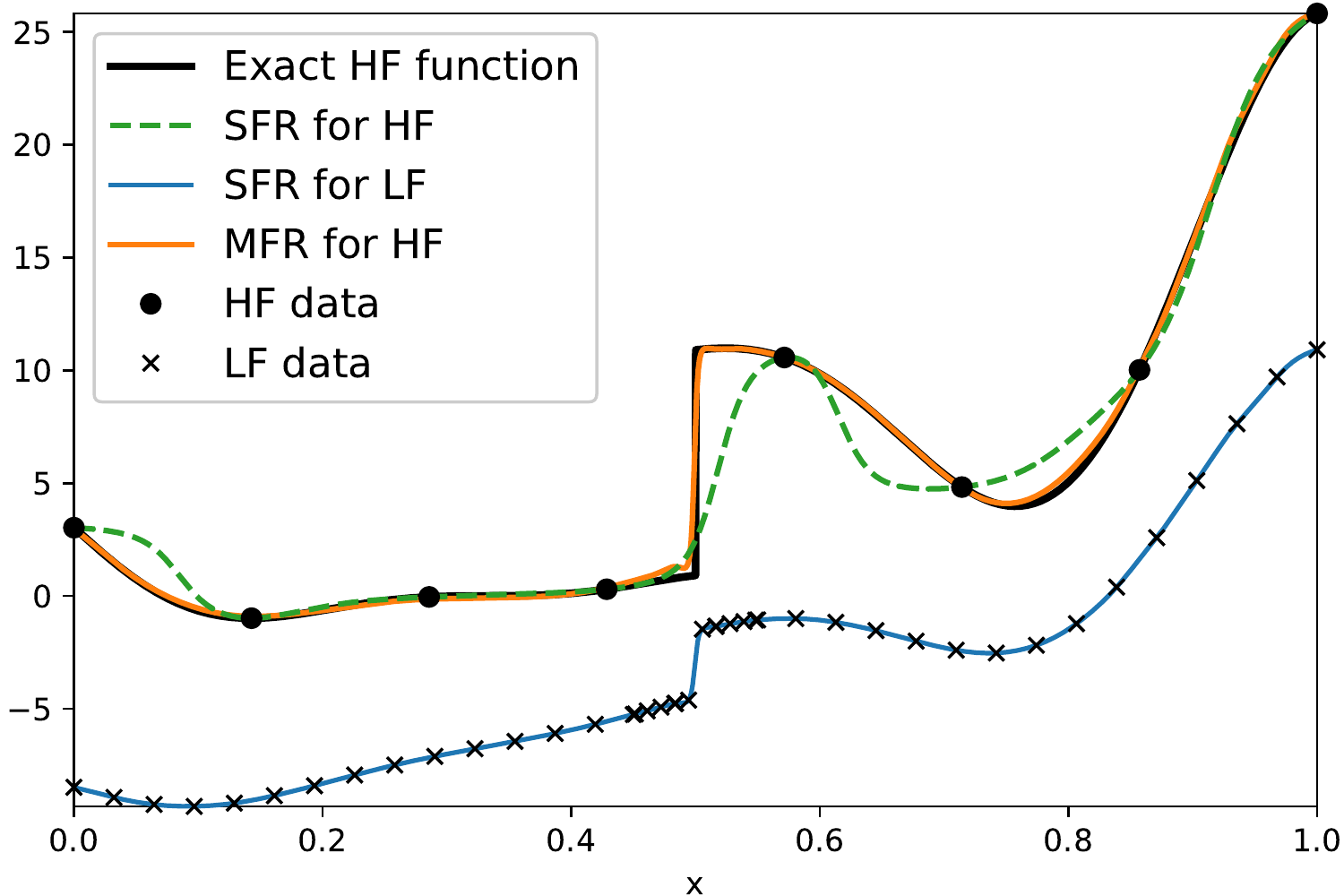}}
	\quad
	\subfigure[MFR using different NNs.]{\includegraphics[width=.475\textwidth]{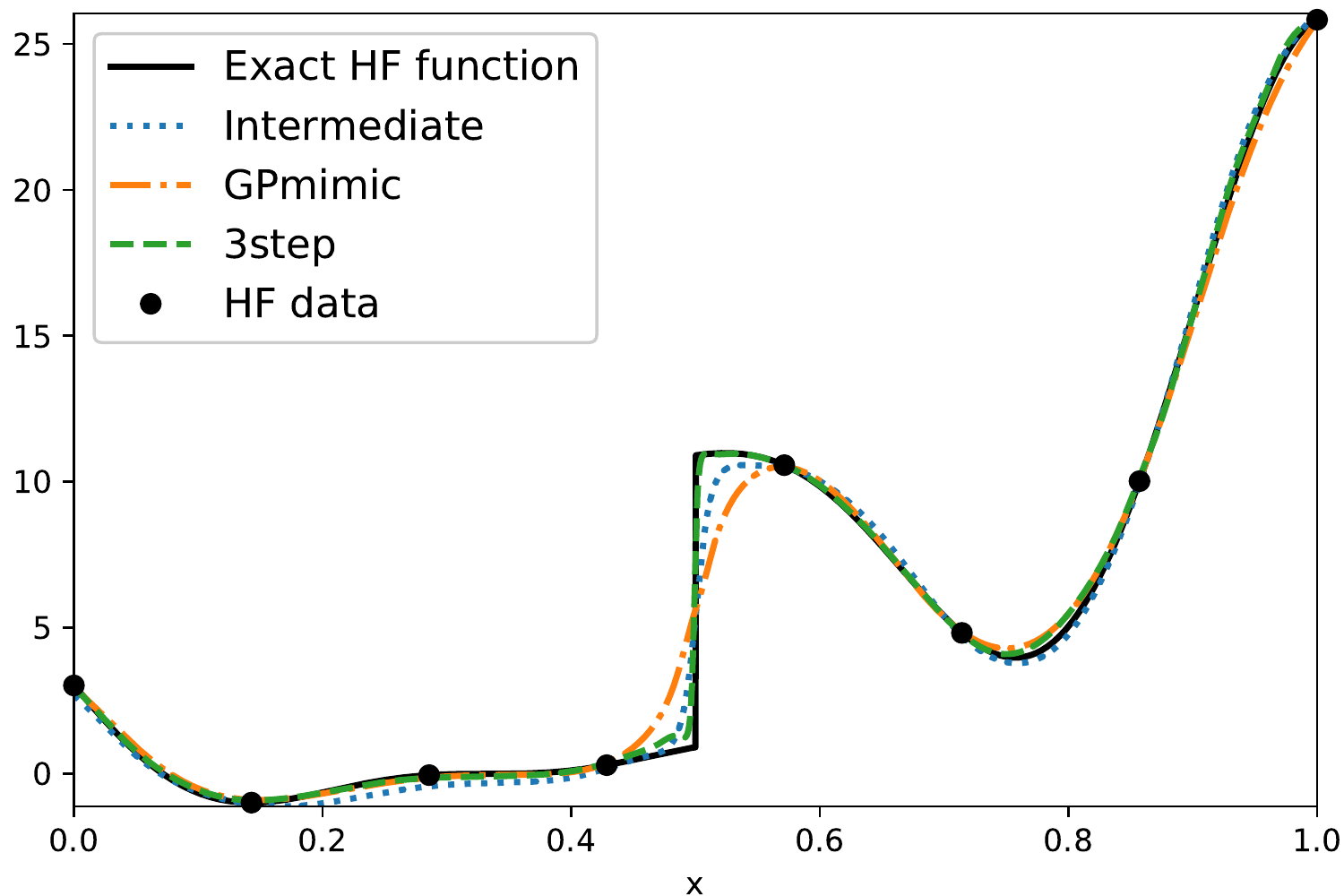}}
	\caption{Discontinuous function: single- and MF regression results using GPs and different ANNs. There are 8 HF (circle) and 42 LF (cross) observations.}
	\label{fig:bench2_mf}
\end{figure}

This second benchmark is designed to analyze how well the proposed NN models can approximate discontinuities in a MF setting. HF and LF data are generated from the following functions:

\begin{align*}
f_\texttt{LF}(x) &= 
\begin{cases}
0.5(6x-2)^{2}\sin(12x-4)+10(x-0.5)-5, & 0 \leq x < 0.5 \\
3+0.5(6x-2)^{2}\sin(12x-4)+10(x-0.5)-5, & 0.5 < x\leq 1 \\
\end{cases}
\\
f_\texttt{HF}(x) &=    
\begin{cases}
2f_\texttt{LF}(x)- 20(x-1), &  0 \leq x < 0.5 \\
4+2f_\texttt{LF}(x)- 20(x-1), & 0.5 < x\leq 1. \\
\end{cases}
\end{align*}
8 and 32 equally spaced locations over $\Omega = [0,1]$ are used as the HF and LF inputs, respectively. In addition, 10 equally spaced points in the interval $[0.45,0.55]$ are added to the LF data set to allow for a better approximation of the discontinuity. In contrast to the first test case, the correlation between the two functions is only piecewise linear, and the piecewise definition results in the discontinuities of distinct amplitudes, see Fig. \ref{fig:bench2_functions}.

Fig.~\ref{fig:bench2_mf} shows that the regression models based merely on the HF data are not able to approximate the discontinuity at $x=0.5$, but the regression results are significantly improved by taking the LF data into account. With $R^2$ scores over 0.99, both the 2-step and 3-step multilevel NN models show a good match with the exact solution. However, the discontinuity is considerably smoothened when other NN models are used, especially when employing the GPmimic model. While the co-kriging presents a discontinuity at $x=0.5$,  fluctuations are present in the interval $[0.4,0.6]$. Despite its low test error, the Intermediate model yields a significantly higher validation error than all the other NNs.  

\begin{table}[h!]
	\caption{Discontinuous function: comparison of the MF regression models. Indicated times account for HPO and final model predictions.}
	\renewcommand{\arraystretch}{1.3}
	\centering%
	\begin{tabular}{c|cccc} 
		\hline
		Model & Validation MSE & Test MSE & $R^2$ & Elapsed time (s)\\
		\hline
		Co-kriging  & -  & 8.07$\times 10^{-1}$  &0.983 & 35\\
		Intermediate & 9.14$\times 10^{0}$  & 2.99$\times 10^{-1}$  & 0.994 & 1587\\
		GPmimic & 1.66$\times 10^{0}$  & 6.92$\times 10^{-1}$  & 0.985 & 1511\\
		2-step & 3.39$\times 10^{1}$  & $\mathbf{9.52\times 10^{-2}}$  & 0.998 & 311\\
		3-step & 7.06$\times 10^{-2}$  & 1.45$\times 10^{-1}$  & 0.997 & 624\\
		\hline
	\end{tabular}
	\label{bench2summary} 
\end{table}

\subsection{Benchmark case 3: Nonlinear correlation}
\begin{figure}[h!]
	\centering
	\subfigure[HF and LF functions.]{\includegraphics[width=.475\textwidth]{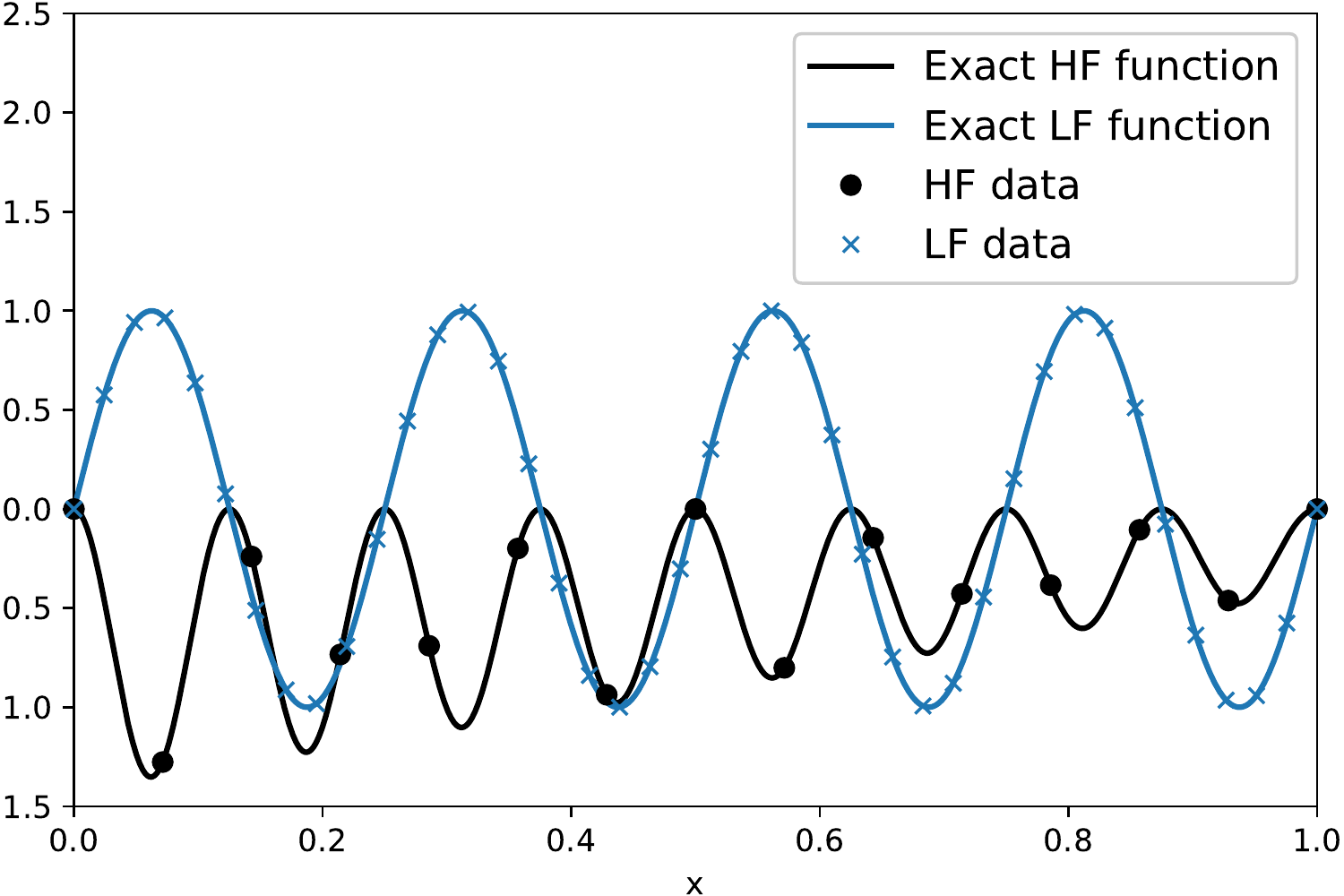}\label{fig:bench3}}
	\quad
	\subfigure[SFR with GP, MFR with co-kriging.]{\includegraphics[width=.475\textwidth]{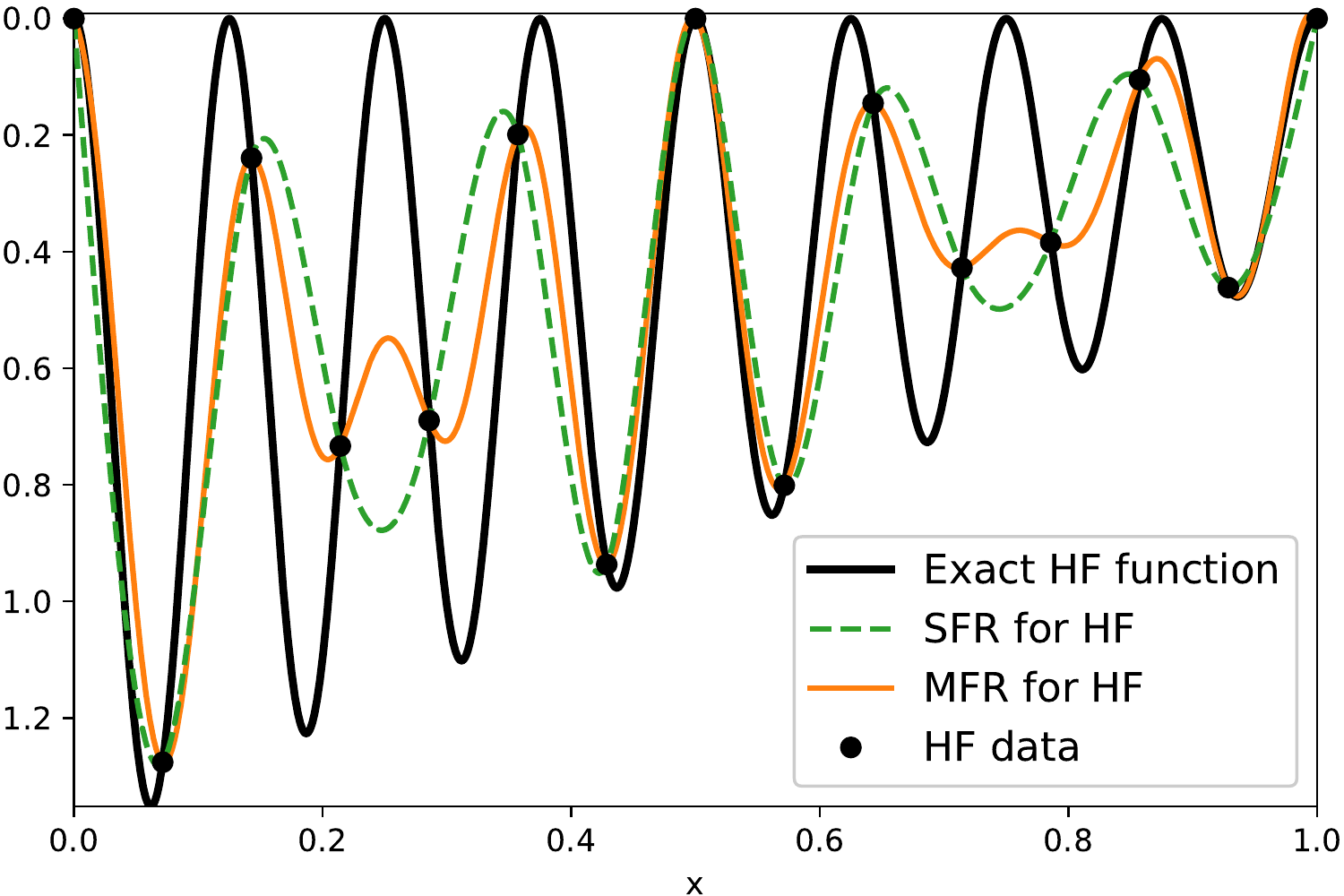}}
	\par\bigskip
	\subfigure[SFR with NN, MFR with 3-step.]{\includegraphics[width=.475\textwidth]{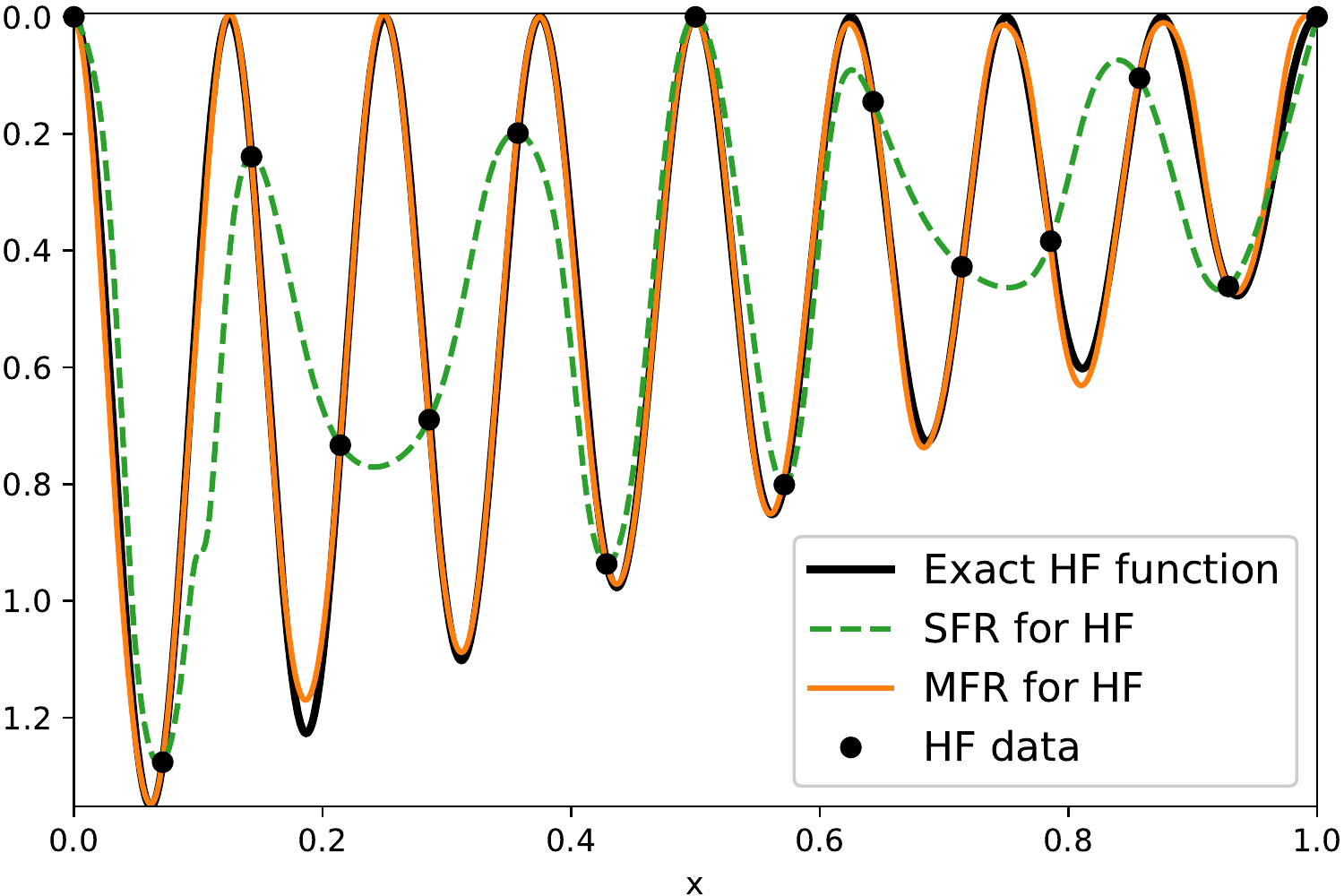}}
	\quad
	\subfigure[MFR using different NNs.]{\includegraphics[width=.475\textwidth]{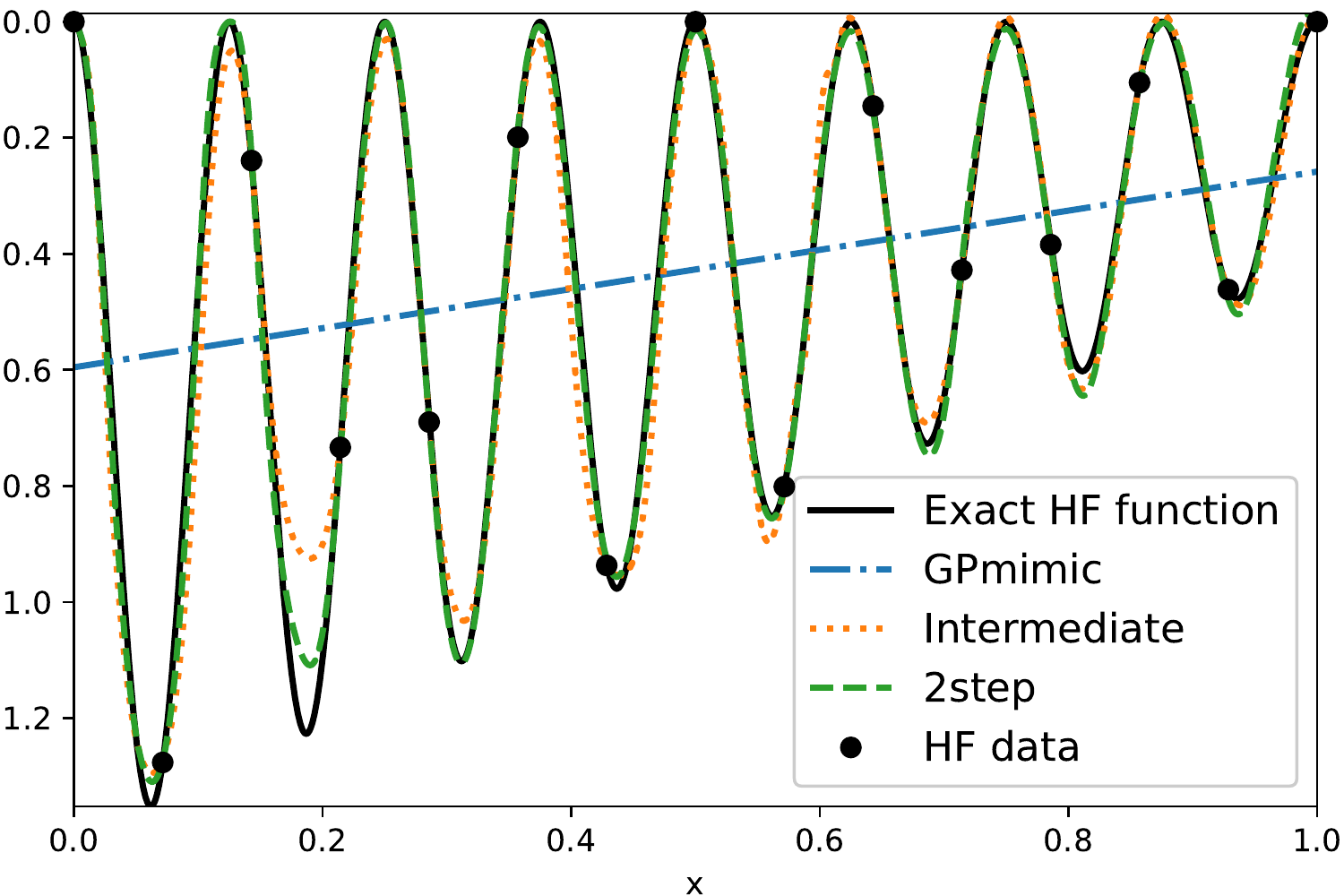}}
	\caption{Nonlinear correlation: single- and MF regression results using GPs  and different ANNs. There are 15 HF (circle) and 42 LF (cross) observations. (a) has a different range for the vertical axis than (b), (c) and (d).}
	\label{fig:bench3_results}
\end{figure}

The third one-dimensional benchmark will test the proposed NN models on a nonlinear correlation between HF and LF levels. The data samples are obtained from the following functions:
\begin{align*}
f_\texttt{LF}(x) &= \sin(8\pi x),\\
f_\texttt{HF}(x) &= (x - \sqrt{2})f_\texttt{LF}^2(x).
\end{align*}
We use 15 equally spaced HF data points and 42 LF data points over the interval $\Omega = [0,1]$. Fig. \ref{fig:bench3} shows that the frequency of the HF function differs from that of the LF function. Moreover, the amplitude of the HF function is linearly decreasing with $x$.

\begin{table}[h!]
	\caption{Nonlinear correlation:  comparison of the MF regression models. Indicated times account for HPO and final model predictions.}
	\renewcommand{\arraystretch}{1.3}
	\centering%
	\begin{tabular}{c|cccc}
		\hline
		Model & Validation MSE & Test MSE & $R^2$ & Elapsed time (s)\\
		\hline
		Co-kriging & -  & $5.94\times10^{-2}$  & 0.561& 36\\
		Intermediate & $8.60\times10^{-2}$  & $5.48\times10^{-3}$  & 0.959& 2135\\
		GPmimic & $1.46\times10^{-1}$  & $1.18\times10^{-1}$  & 0.128 & 2019\\
		2-step & $9.51\times10^{-4}$  & $8.09\times10^{-4}$  & 0.994& 1072\\
		3-step & $1.38\times10^{-3}$ & $\mathbf{4.49\times10^{-4}}$  & 0.997 & 1208\\
		\hline
	\end{tabular}
	\label{bench3summary}
\end{table}

As in the previous test cases, single-fidelity regression models do not provide satisfactory accuracy. Neither co-kriging nor the GPmimic model can leverage the LF data to improve the regression results. With an $R^2$ score of 0.128, the GPmimic NN hardly performs better than the average function value. However, it is not surprising that both models fail in the current benchmark, as their mathematical structure is unable to detect nonlinear correlations between the two fidelity levels. Fig. \ref{fig:bench3_results} underlines that the difficulty of this test case lies in approximating the different local extrema of the HF function, since data points are not always available close to the local extremum. Nevertheless, both the multilevel NN models recover an $R^2$ score exceeding 0.99. 

\subsection{Benchmark case 4: 20-D benchmark}

It is well known that GPR suffer from the curse of dimensionality, and cannot be effectively used when the number of data points is very large $(N > 10000)$, especially when big data sets are needed to sufficiently cover the high-dimensional input space. The current benchmark is chosen to show that, in contrast to GPR, NNs remain a valid candidate for MF regression in the presence of high dimensionality and large data sets. 
The following 20-dimensional functions \cite{Meng} define the MF setting in this example:
\begin{align*}
f_\texttt{HF}(\mathbf{x}) &= (x_1 - 1)^2 + \sum_{i=2}^{20}(2x_i^2 - x_{i-1})^2\,,\\
f_\texttt{LF}(\mathbf{x}) &= 0.8f_\texttt{HF}(\mathbf{x}) - \sum_{i=2}^{20}0.4x_{i-1}x_{i} -50\,,
\end{align*}
with $\vb{x}=\{x_1,x_2,\cdots,x_{20}\}\in \Omega = [-3,3]^{20}$. HF  data are sampled from $f_\texttt{HF}$ at 5000 locations randomly chosen from a uniform distribution over $\Omega$, while the LF  samples are evaluated at 30000 random input locations.

\begin{figure}[h!]
	\centering
	\subfigure[Single-fidelity regression]{	\includegraphics[width=.425\textwidth]{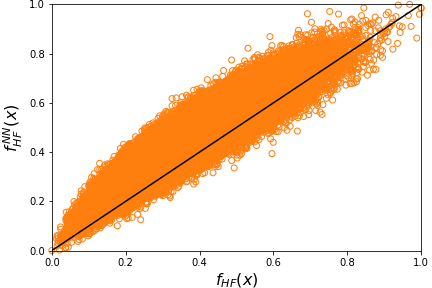}\label{fig:20Dsf}}
	\quad
	\subfigure[MF regression with the 3-step model]{\includegraphics[width=.425\textwidth]{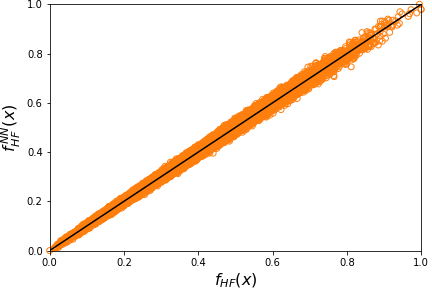}\label{fig:20Dmf}}
	\caption{20-D benchmark: Exact function values versus predicted values by NNs at one million random locations. All data are normalized to the range [0,1].}
	\label{fig:20D}
\end{figure}

In this case, the co-kriging would require one to compute the inverse of a $35000\times35000$ matrix. Just the storage of such a matrix can take around 10 GB of memory. Such extensive computational cost makes it impossible use a GPR model. Instead we utilize the proposed 3-step NN model whose training is accelerated by a GPU, and we opt for the 5-fold cross-validation to tune NN hyperparameters. Fig. \ref{fig:20D} shows that the 3-step model can accurately predict the value of $y_\texttt{HF}$ at one million random input locations. 

\section{Numerical results (II): application to parametrized PDEs}

In this section we apply our MF framework to a problem arising in acoustics, namely the propagation of a pressure wave $P(\textbf{x},t)$ into an acoustic horn with parametrized shape, addressed in \cite{Horn_Manzoni}. In particular, we are interested in evaluating an input-output map involving the solution of a PDE, exploiting the accuracy of the finite element method to obtain HF solutions, while relying on a reduced basis (RB) method to compute LF, but fast and inexpensive, approximations. Different from the previous section, we consider different sources of data for each fidelity level. 

\begin{figure}[h]
    \centering
    \includegraphics[width=0.49\textwidth]{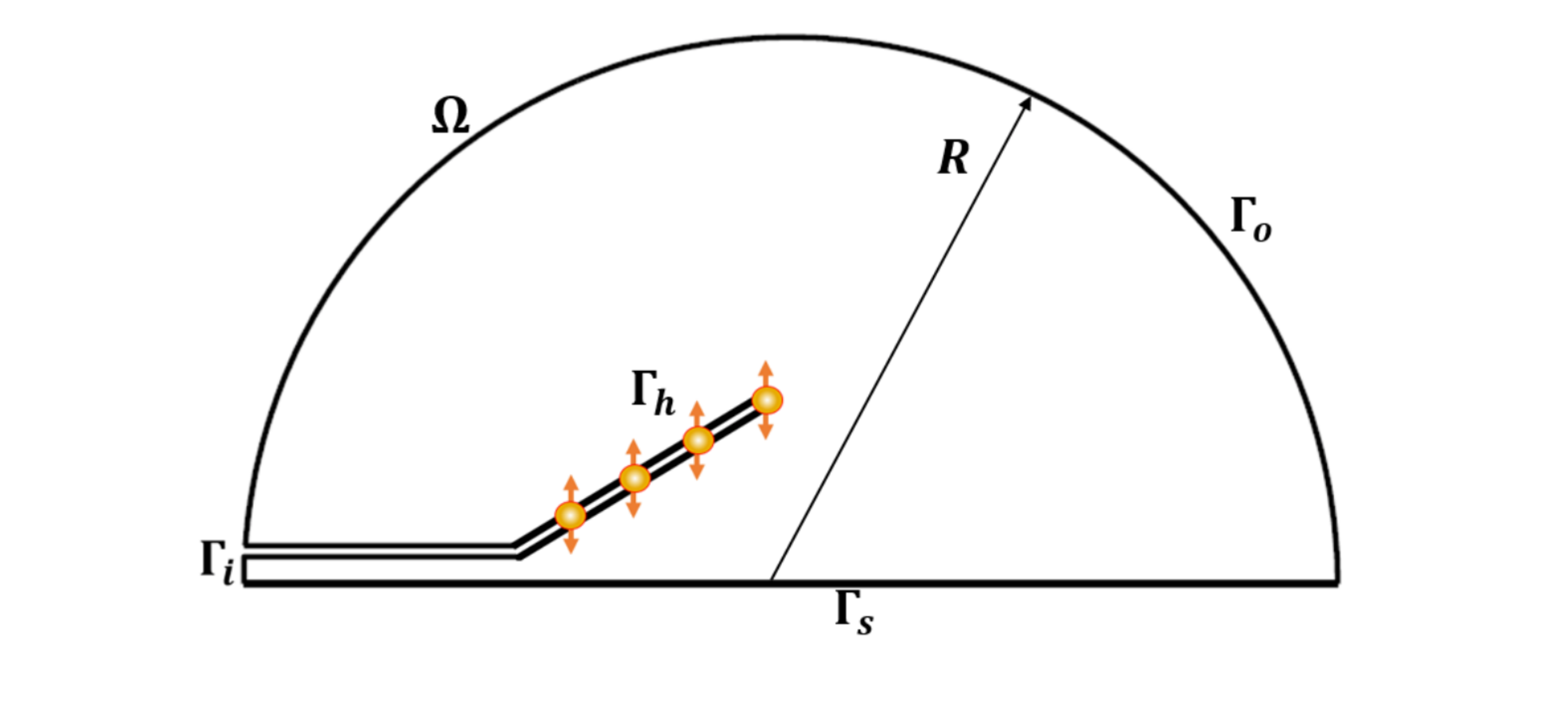}
    \includegraphics[width=0.49\textwidth]{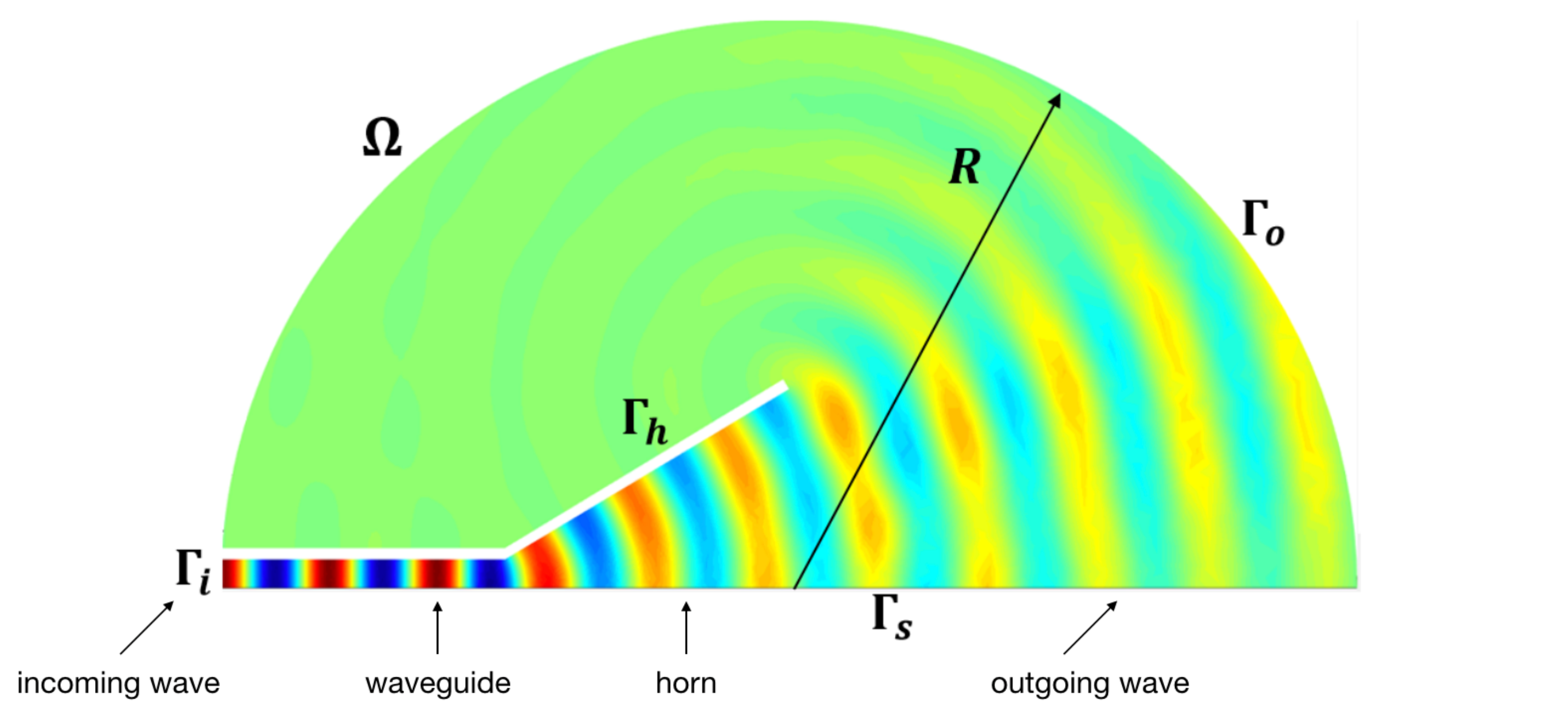}
    \caption{Acoustic horn problem. Left: the computational domain $\Omega$ and the boundaries. Right: the control points used in RBF shape parametrization, whose vertical displacements are treated as parameters.}
    \label{fig:geom_horn}
\end{figure}

We consider an acoustic device, illustrated in Fig \ref{fig:geom_horn} (left), comprising of a waveguide, with infinite extension to the left and a conical extremity on the right, namely the horn\footnote{For simplicity, the device extends infinitely in the direction normal to the plan and its wall consists of sound-hard material. Therefore, for the frequencies in the range under consideration, we assume that all non-planar modes in the waveguide are negligible, which allows us to reduce the problem to two space dimensions.}, and an internal propagating planar wave inside the waveguide: once the wave reaches the horn, a portion of its energy is converted into an outer-going wave. Under the assumptions of single-frequency and time harmonic waves, the acoustic pressure can be expressed as $P(\textbf{x},t) = \mathfrak{R}(p(\textbf{x})e^{i\omega t})$, where the complex amplitude $p(\textbf{x})$ satisfies the monochromatic steady-state Helmholtz equation with mixed Neumann-Robin boundary conditions: 
\begin{alignat}{2}
    \Delta p + k^2p &= 0 \qquad  &&\text{in } \Omega \nonumber \\
    (ik + \frac{1}{2R})p + \nabla p \cdot \textbf{n} &= 0  \qquad &&\text{on }\Gamma_{0} \nonumber\\
    ikp + \nabla p \cdot \textbf{n} &= 2ikA \qquad &&\text{on } \Gamma_{i} \nonumber  \\
    \nabla p \cdot \textbf{n} &= 0  \qquad &&\text{on } \Gamma_{h} \cup \Gamma_{s} = \Gamma_{n},  
    \label{eq: Helmholtz_system}
\end{alignat}
where $k = \omega / c$ is the wave number, $\omega = 2\pi f$ the angular frequency and $ c = 340 \text{ cm s}^{-1}$ the speed of sound;  $\textbf{n}$ denotes the outward-directed unit normal on the boundary of $\Omega$. We restrict the computation to the domain $\Omega$ shown in Fig.~ \ref{fig:geom_horn}, and impose on $\Gamma_{i}$ -- a propagating wave with amplitude $A=1$ while absorbing the outer-going planar waves -- an absorbing condition on the far-field boundary $\Gamma_{o}$, and  homogeneous Neumann boundary conditions on the sound-hard walls of the device $\Gamma_{h}$ as well as on the symmetry boundary $\Gamma_{s}$. We take, for simplicity, the radius equal to $R=1$ \cite{bangtsson2003shape}.

In addition to the frequency $f$ of the incoming wave, we parametrize, as in \cite{Horn_Manzoni}, the shape of the horn by means of radial basis functions, introducing as parameters $\boldsymbol{\mu}_g = (\mu_{g,1}, \ldots, \mu_{g,4})$ the vertical displacement of four control points located on the horn wall $\Gamma_{h}$, shown in Fig. \ref{fig:geom_horn}, right. The admissible domain configurations are defined as the diffeomorphic images $\Omega(\boldsymbol{\mu}_g)$ of the reference shape $\Omega$ through a deformation mapping  $\mathbf{T}(\cdot;\boldsymbol{\mu}_g)$ obtained as linear combinations of the control points displacements. Hence, the acoustic problem depends on five parameters, i.e., we denote by $\boldsymbol{\mu} = ( f, \boldsymbol{\mu}_{g} )$ the parameter vector and  by $\mathcal{D} \subset \mathbb{R}^p$ the parameter space.

We focus our analysis on a specific output of interest, namely the index of reflection intensity (IRI) \cite{bangtsson2003shape} which measures the transmission efficiency of the acoustic horn and is defined as the absolute value of the average reflected wave at the sound inlet $\Gamma_i$, i.e.,
\begin{equation*}
     f(\boldsymbol{\mu}) = J(p(\boldsymbol{\mu})) =  \left| \frac{1}{|\Gamma_i|} \int_{\Gamma_i}p(\boldsymbol{\mu})\,d\Gamma -1 \right|. 
\end{equation*}
In particular, 
\[
f_\texttt{LF}^N(\boldsymbol{\mu}) = 
J(p_N^m(\boldsymbol{\mu})), \qquad  \mbox{and} \qquad
f_\texttt{HF}(\boldsymbol{\mu}) = 
J(p_h(\boldsymbol{\mu})), 
\]
define the MF setting for this test case. We denote by  $p_N^m(\boldsymbol{\mu}))$  and $p_h(\boldsymbol{\mu}))$ the LF and the HF solutions of the parametrized PDE, respectively, the construction of which is reported in the Appendix B. Moreover, we highlight the dependence of the LF model on the dimension $N$ of the ROM used to evaluate the PDE solution. In the following subsections, we consider two different scenarios, dealing with either 1 or 5 parameters.

\subsection{Case with $p=1$ parameter} 

In this first case, we compare different MF strategies, by considering the output of interest as a function of the frequency $\boldsymbol{\mu} = f$ only, letting it vary in $\mathcal{D} = [10,1800]$. Hence, we first limit our analysis to the reference configuration of the horn, without taking geometric parameters into account. 
In this specific case, the linear system arising from the FE approximation of the problem \eqref{eq: Helmholtz_system} exhibits an affine decomposition (see Appendix B), so that any further hyper-reduction stage is not required when constructing the ROM. To build this, we first randomly sample 150 values of the frequency $f$ and compute the corresponding HF solutions through the FOM. The FOM is approximated by $\mathbb{P}_1$ finite elements and, considering a mesh made of 8740 triangular elements, we have an HF model of dimension $n = 4567$. Then, we apply POD and extract $N=44$ reduced basis functions by imposing a relative projection error of $10^{-5}$. All computational details are summarized in Table \ref{tab: Comp_det_1}. The computation of HF and LF solutions, i.e. the FOM and ROM solutions, respectively, is carried out in \textsc{Matlab}, using the redbKIT library \cite{redbKIT}. All hyperparameters obtained by HPO are reported in Appendix A.

\begin{figure}[t!] 
    \centering
    \includegraphics[width=.8\textwidth]{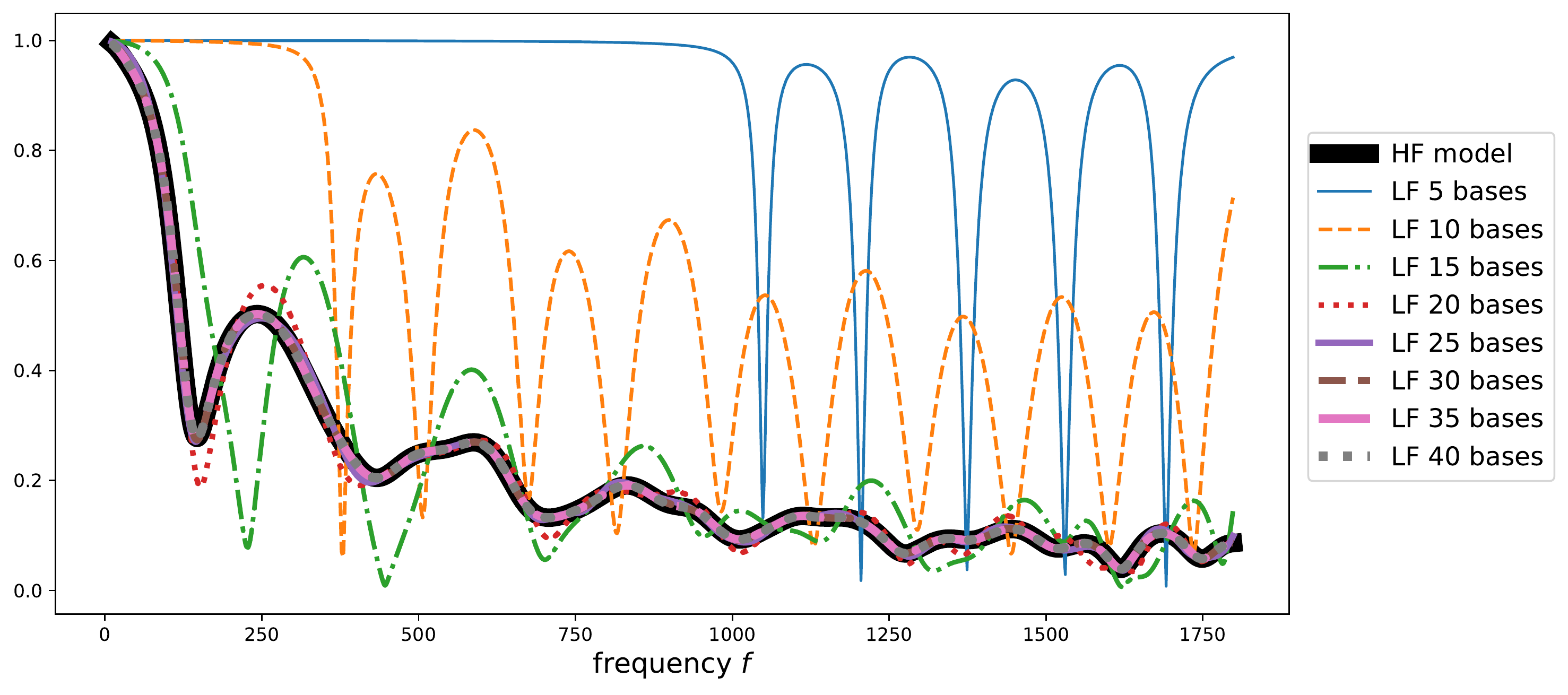}
    \caption{HF model $f_\texttt{HF}(\boldsymbol{\mu})$ and different LF models $f_\texttt{HF}(\boldsymbol{\mu})$ depending on the number $N$ of basis functions. For too small values of $N$, the LF model prediction is rather poor, while it is almost indistinguishable from the HF model prediction for $N \geq 25$.}
    \label{fig: LF_quality}
\end{figure}

We assess how the quality of the LF model and the amount of HF data impact the accuracy of the MF prediction. We recall that {\em (i)} we can improve the quality of the LF model by selecting a smaller or larger number of bases, and  {\em (ii)} we can freely decide the parameter values for which we solve the HF model, by keeping the sampling method fixed. A comparison among LF models obtained with different dimensions $N$ of the ROM is reported in Fig.~\ref{fig: LF_quality}. We choose to limit the selected basis functions in the POD-Galerkin ROM to be between 5 and 22 (even though the maximum number of available basis functions is 44), to deal with a potentially inaccurate (or, at least, not accurate enough) LF model. The ability to obtain accurate predictions by combining a few HF data and several evaluations of a reliable, but not sufficiently accurate, ROM, is indeed an attractive feature of the proposed framework, as this may prevent us from constructing ROMs with large dimensions and poor efficiency. Moreover, we consider the same amount of HF data, while we keep fixed (and equal to 32) the number of LF data.

\begin{table}[h]
    \centering
    \caption{Computational details in the case with $p=1$ parameter: $\boldsymbol{\mu} = f$.}
    \vspace{5pt}
    \begin{tabular}{l*{6}{l}r}
    \hline
    Number of parameters            & 1 & Parameter domain $\mathcal{D}$ & $[10,1800]$   \\
    Number of FE elements          & 8740 & Tolerance RB POD & $10^{-5}$  \\
    Number of FE dofs n & 4567 & Number of ROM dofs N         & 44   \\
    Number of HF data & from 5 to 22 & Number of bases & from 5 to 22 \\
    Number of LF data & 32 & Sampling method & LHS\\
    \hline
    \end{tabular} 
    \label{tab: Comp_det_1}
\end{table}

\begin{figure}[t!]
    \centerline{
    \includegraphics[width=.9\textwidth]{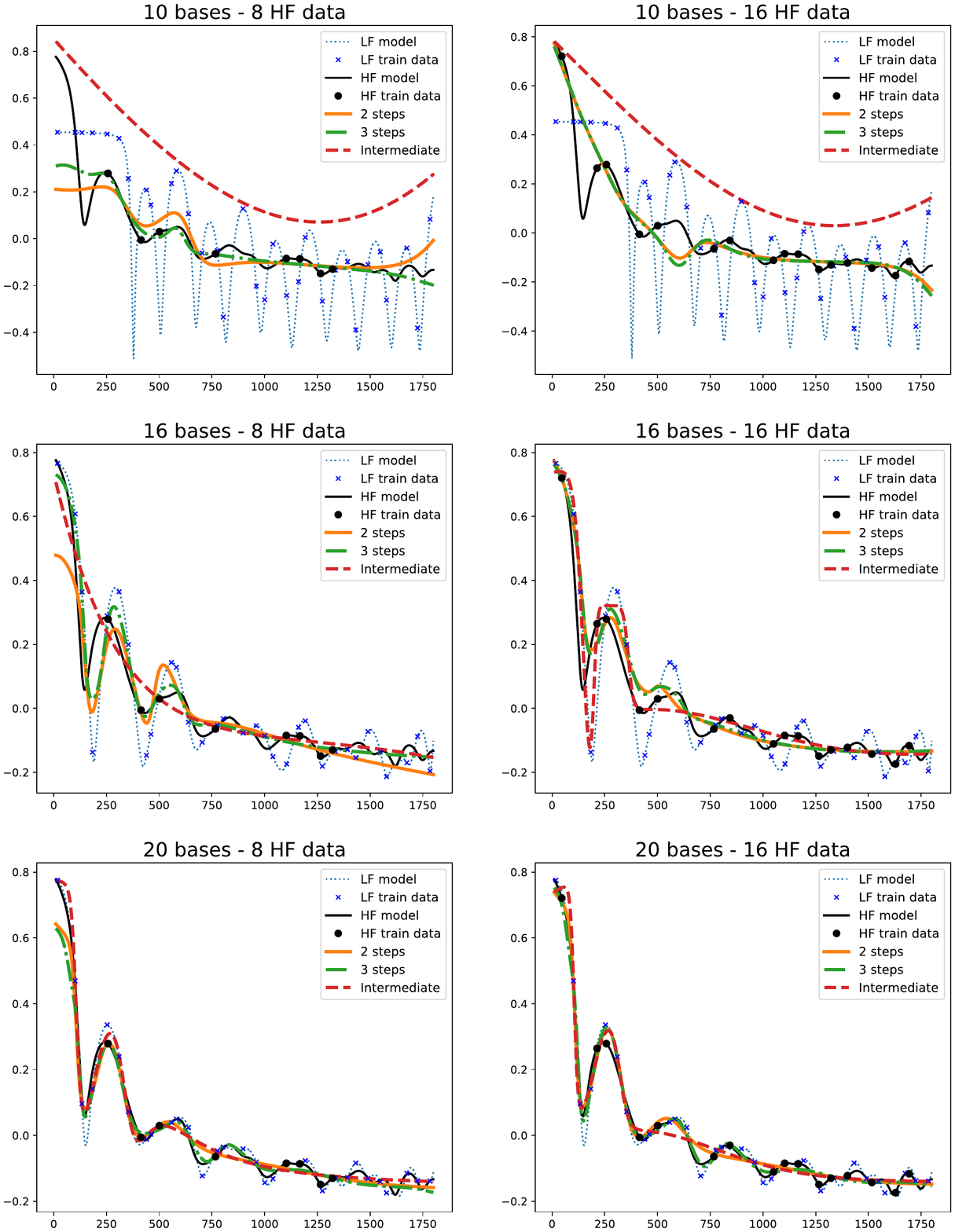}
    }
    \caption{Plots of the HF, LF and MF (\textit{Intermediate, 2-step, 3-step}) models considering different number of bases and HF training data.}
    \label{fig: graph_1param}
\end{figure}

\begin{figure}[t!]
    \centerline{
    \includegraphics[width=\textwidth]{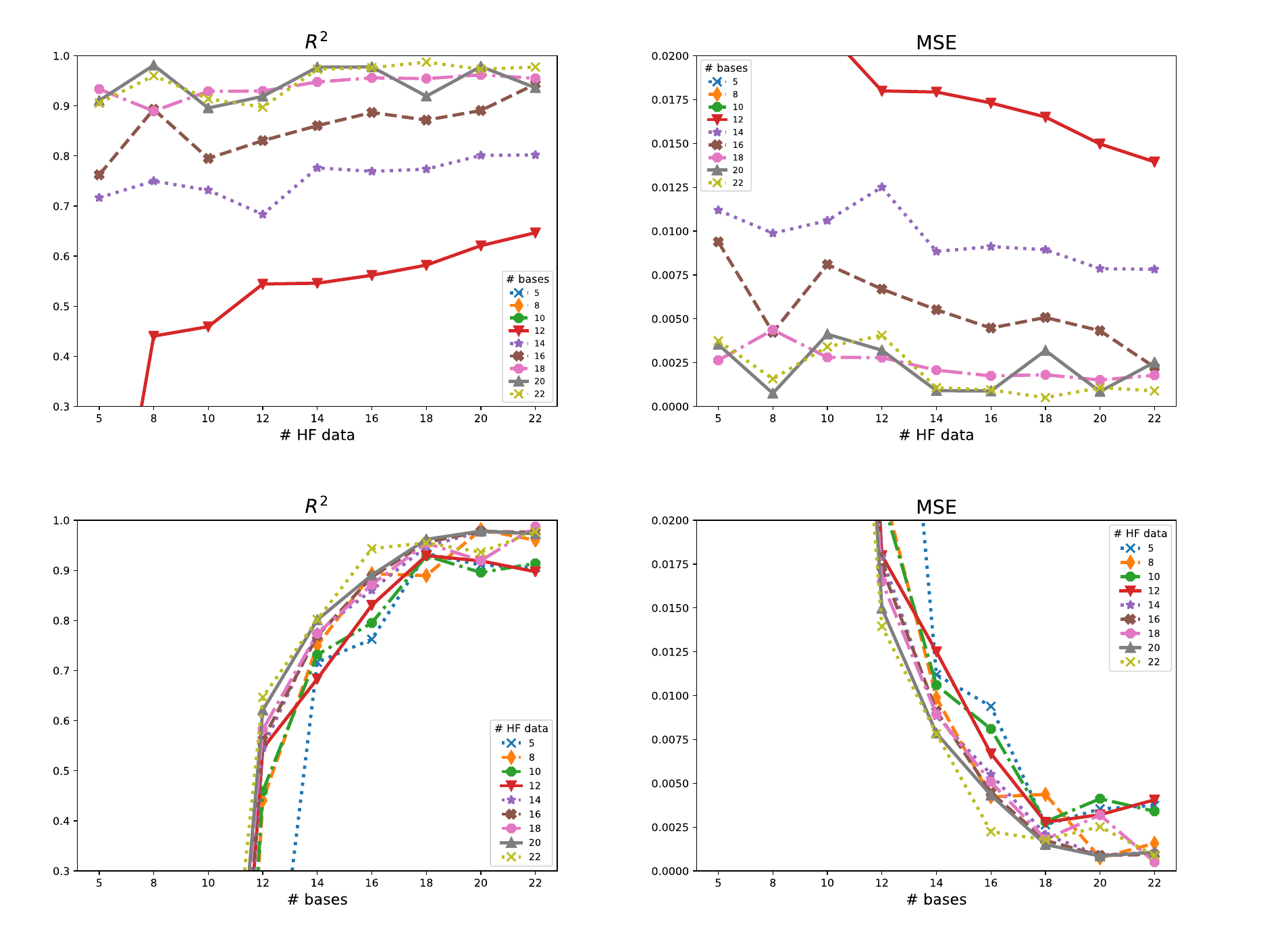}
    }
    \caption{Case with $p=1$ parameter, \emph{Intermediate} architecture. MSE and $R^2$ for different amounts of HF data and LF dimension $N$.}
    \label{fig: r2_mse_1param_int}
\end{figure}

\begin{figure}[t!]
    \centerline{
    \includegraphics[width=\textwidth]{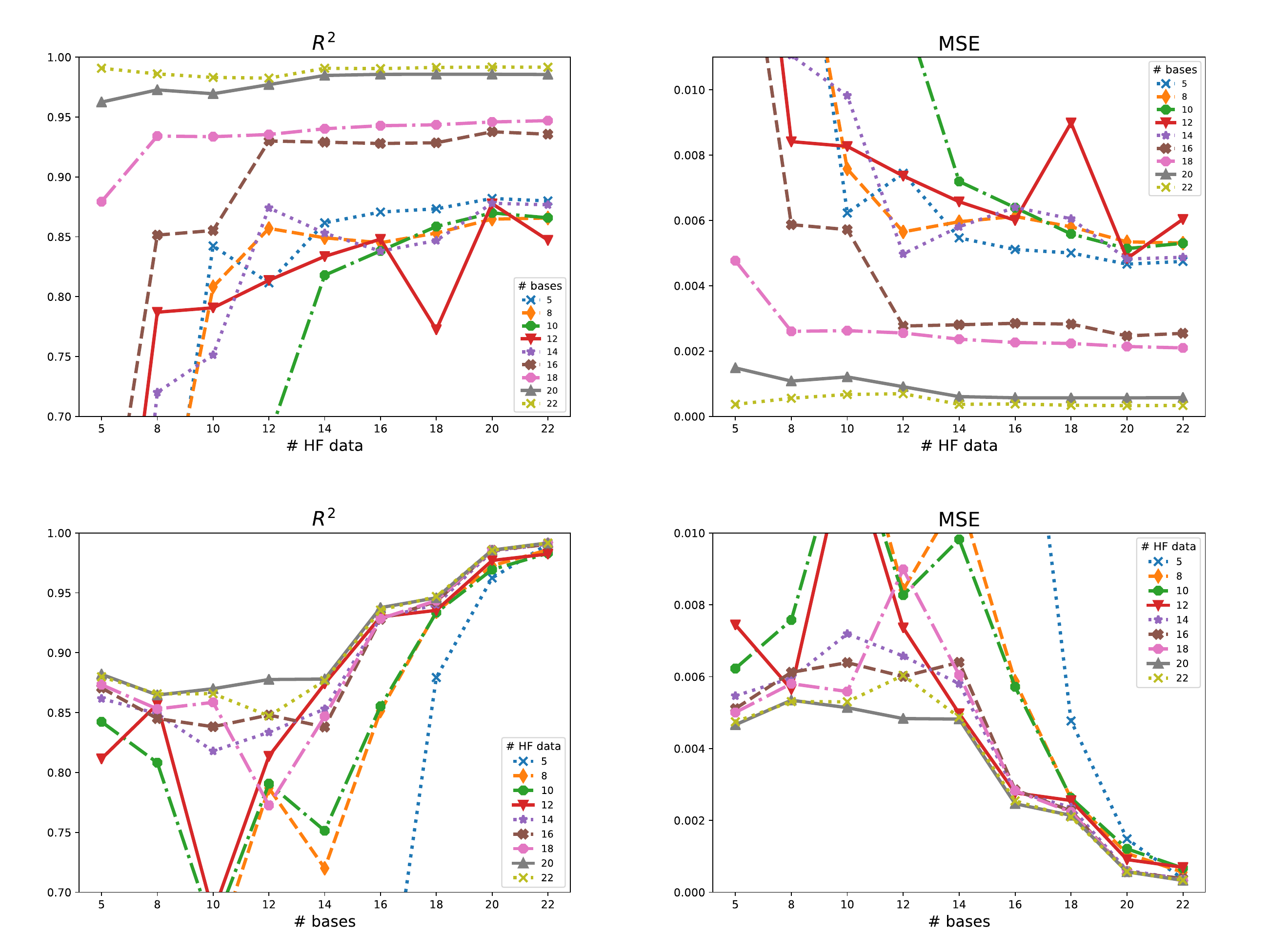}
    }
    \caption{Case with $p=1$ parameter, \emph{2-step} architecture. MSE and $R^2$ for different amounts of HF data and LF dimension $N$.}
    \label{fig: r2_mse_1param_2s}
\end{figure}

\begin{figure}[t!]
    \centerline{
    \includegraphics[width=\textwidth]{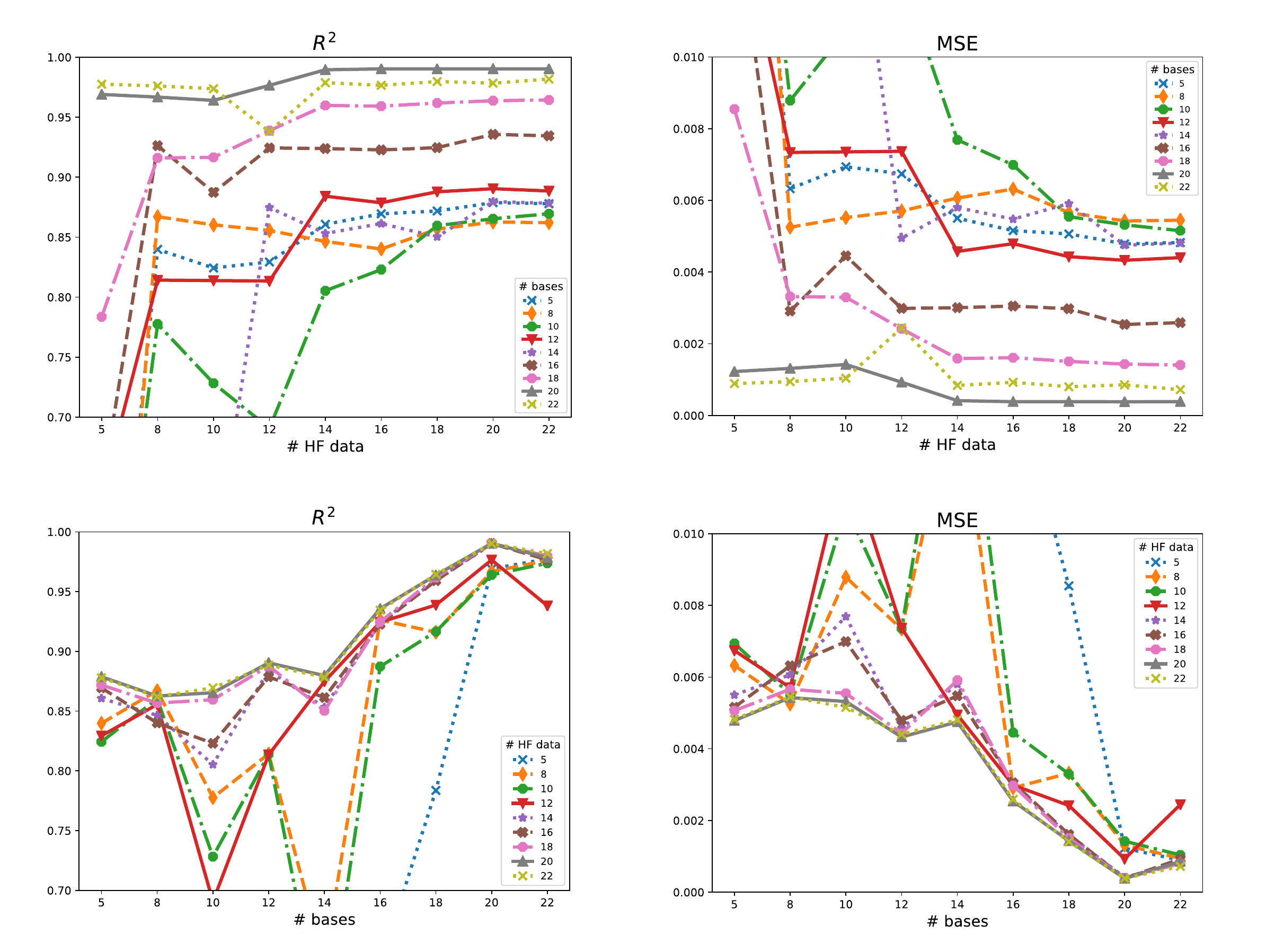}
    }
    \caption{Case with $p=1$ parameter, \emph{3-step} architecture. MSE and $R^2$ for different amounts of HF data and LF dimension $N$.}
    \label{fig: r2_mse_1param_3s}
\end{figure}

In the following, we test the \textit{Intermediate, 2-step, 3-step} network architectures; for each combination \textit{number of HF data - number of bases}, we train the NNs, predict the output of interest and compute the indices of goodness of fit $R^2$ and $MSE$. Evaluated outputs with these network architectures are displayed in Fig.~\ref{fig: graph_1param}; the values of $R^2$ and $MSE$ are reported instead in Figs.~\ref{fig: r2_mse_1param_int}, \ref{fig: r2_mse_1param_2s} and \ref{fig: r2_mse_1param_3s},  as functions of the number $N$ of basis functions, and the amount of HF data considered.

Overall, the neural architectures perform very well and provide good predictions, in terms of both $MSE$ and $R^2$, provided a sufficient number of basis functions and HF data are used. 
The \textit{2-step} and \textit{3-step} models produce similar results and perform better than the \textit{Intermediate} model, as displayed in Fig.~\ref{fig: graph_1param}. In particular, the multilevel models are robust and efficient even in the cases where the LF models are built from few bases, while the \textit{Intermediate} network has poor predictive accuracy without a sufficiently accurate LF model, even if a large amount of HF data are provided. In fact, with multistep networks we can reach values of $R^2$ larger than 0.8 even considering just $N=5$ bases, whereas the \textit{Intermediate} model does not provide large values of $R^2$ with less than 12 bases, see Figs.~\ref{fig: r2_mse_1param_int}, \ref{fig: r2_mse_1param_2s} and \ref{fig: r2_mse_1param_3s}. Regarding the computed outputs, we see how the peaks with larger amplitude found for $f <1000$ are correctly described by  both the   \textit{2-step} and \textit{3-step} models, while smaller amplitude peaks for $f>1000$ are better captured by the \textit{3-step} model than the \textit{2-step}. On the other hand, the intermediate model provides a less accurate trend of the output, and only provides reliable results provided that both $N$ and the amount of HF data are large enough. 

As expected, the prediction improves as the number of HF data or the number of basis functions increases, even if these two features impact the accuracy in a slightly different way. Restricting, for the sake of simplicity, to the case of a \emph{3-step} model (see Fig.~\ref{fig: r2_mse_1param_3s}, top) we observe that,
for each fixed number of reduced bases $N$, the goodness-of-fit indices improve as the number of HF data increases until a threshold limit is reached, which is determined solely by the number of bases and cannot be overcome by adding more HF data.
Conversely (see Fig.~\ref{fig: r2_mse_1param_3s}, bottom) by increasing  $N$, the trends of $R^2$ (resp. MSE), corresponding to the different numbers of HF data adopted, increase (resp. decrease) and also tend to reach values closer to each other. The amount of HF data thus becomes less and less important as the LF model improves. Hence, in this specific problem it is more efficient to have a good LF model even with a small amount of HF data, rather than lots of HF data but a poor LF model.

 \subsection{Case with $p=5$ parameters} 

We finally apply our MF setting to the case where all the five parameters $\boldsymbol{\mu} = (f,\boldsymbol{\mu}_g)$, namely the frequency and the four geometric parameters, vary, in order to consider shape variations in the horn geometry as well. The parameter space is $\mathcal{D}=[50,1000] \times \mathcal{D}_g$, where $\mathcal{D}_g = [-0.03,0.03]^4$.

We compute 200 FOM snapshots for 200 points sampled in the parameter domain $\mathcal{D}$ through a Latin Hypercube sampling design. In this case, the reduced basis is made from $N = 80$ POD modes; see Table~\ref{tab: Comp_det_5}. As in the case of $p=1$ parameter, we select an appropriate range for the number of bases and HF data to assess the efficiency of the MF approach in different scenarios. We consider a larger number of LF data (500 instead of 100) than in the case of $p=1$ parameter; we then vary the number of HF data from 5 to 45, and the number of basis functions of the LF model from 5 to 40. We report the results only in the case of the \textit{3-step} model for the sake of brevity and to focus on the role of the quality and quantity of the training data rather than on the choice of the NN architecture.

\begin{table}[h]
    \centering
    \caption{Computational details in the case with $p=5$ parameters: $\boldsymbol{\mu} = (f, \boldsymbol{\mu}_g)$.}
    \vspace{5pt}
    \begin{tabular}{l*{6}{l}r}
\hline
    Number of parameters            & 5 & Parameter domain $\mathcal{D}$ & $[50,1000] \times [-0.03,0.03]^4$  \\
    Number of FE elements          & 8740 & Number of FE snapshots & 200  \\
    Number of FE dofs n & 4567 & Number of ROM dofs N         & 80   \\
    Number of HF data & from 5 to 45 & Number of bases & from 5 to 40 \\
    Number of LF data & 500 & Sampling method & LHS\\
    \hline
    \end{tabular} 
    \label{tab: Comp_det_5}
\end{table}

\begin{figure}[h!]
    \centerline{
    \includegraphics[width=\textwidth]{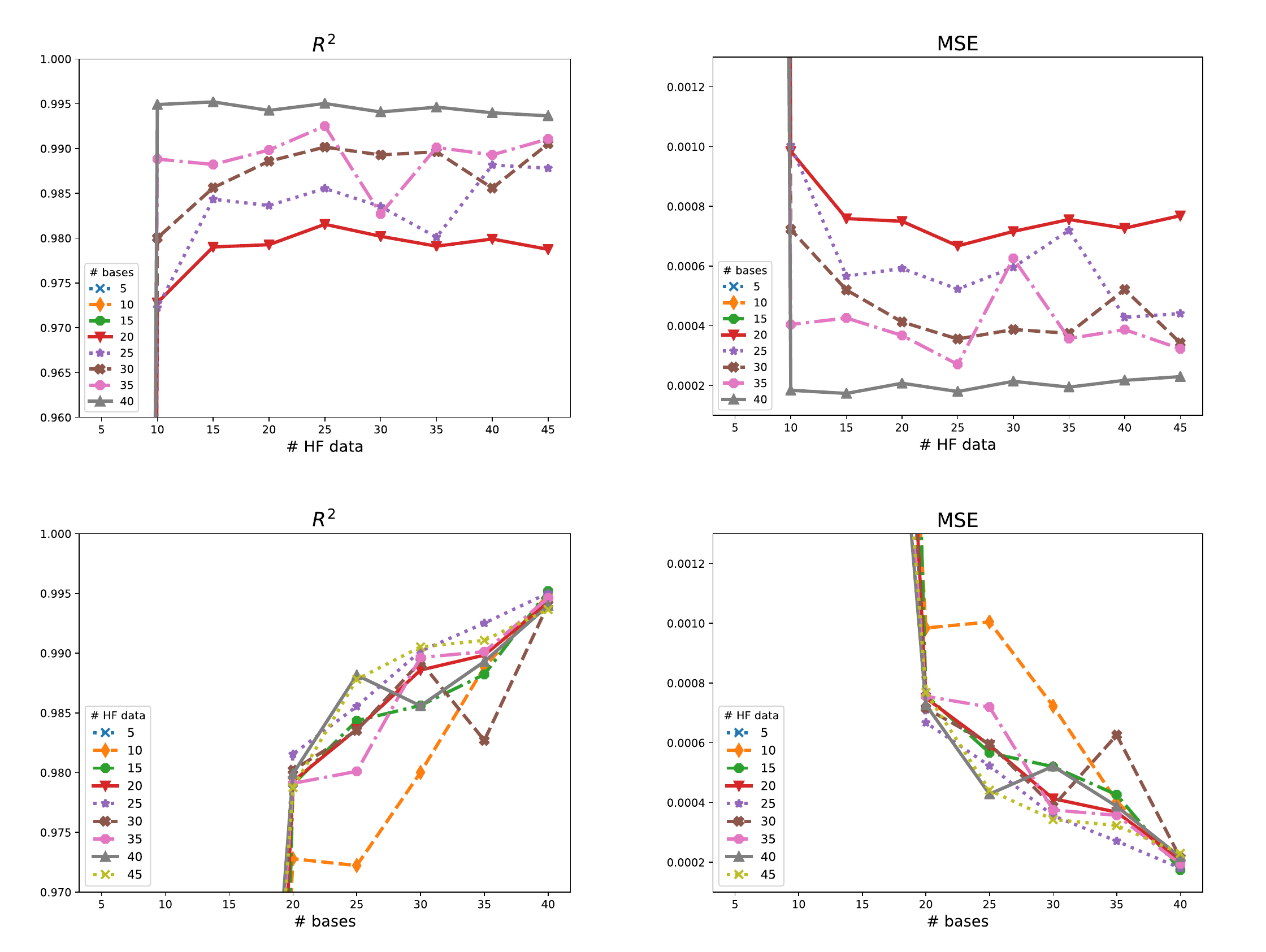}
    }
    \caption{Case with $p=5$ parameters, \emph{3-step} architecture. MSE and $R^2$ for different amounts of HF data and LF dimension $N$.}
    \label{fig: r2_mse_5param_1}
\end{figure}

Passing from 1 to 5 parameters does not worsen the prediction power of the NN architecture. Indeed, once again we obtain an accurate prediction from a small number of FOM data, by exploiting a large number of ROM solutions that can be computed very quickly and inexpensively. As in the previous case, the MF prediction improves both as the number of HF data and the number of bases increase, with the latter playing a more important role in this framework. Similar to  the case $p=1$, we display $R^2$ and the MSE obtained with the \textit{3-step} NN architecture as functions of the aforementioned factors (see Fig.\ref{fig: r2_mse_5param_1}).

Numerical results show that we can achieve very good results even in 5 dimensions by employing 500 LF data, and that increasing both the amount of HF data and the dimension of the LF model improves the prediction accuracy. 
When considering a poor LF model, the number of HF data is fundamental to obtain good prediction: for instance, with an LF  model of dimension $N=15$, we obtain  $R^2 = 0.12768$ with  5 HF data, and $R^2 =0.91093$ with 40 HF data. In contrast, the amount of HF data loses importance as the LF model becomes more accurate. In fact, with an LF model of dimension $N=35$, when HF data increases from 10 to 45, $R^2$ only improves by $0.2\%$, passing from 0.98882 to 0.99107. For a fixed, small number of basis functions (e.g., $N=5$) $R^2$ (resp., the MSE) continues to increase (resp., decrease) as the number of HF data increases, while both indices flatten when considering larger values of $N$ (e.g., $N>15$). Considering a fixed number of HF data greater than 5, it is possible to reach excellent values of the goodness-of-fit indices just by increasing the number of bases; in particular, starting from 20 bases, we already reach a $R^2$ greater than 0.97 and an MSE smaller than $1.2\times10^{-3}$, with both indices continuing to improve as the number of bases increases. In particular, it seems that increasing the number of HF data only works provided also the dimension $N$ of the LF model is enlarged. On the other hand, as $N$ increases, the indices keeps improving and the number of HF data becomes less and less relevant.
Therefore, also in the case of $p=5$ parameters we can conclude that improving the quality of the LF model is a more efficient strategy - as well as computationally cheaper -  than increasing the number of HF data, to reach a certain degree of accuracy.

\section{Conclusions}

In this work we discuss MF regression with ANNs, for which four different architectures are presented. The proposed NN schemes are benchmarked against co-kriging on a collection of test cases of increasing complexity. We also successfully predict output quantities associated with parametrized PDEs using the multi-fidelity NN models. Observations made from different numerical examples show that MF models based on NNs can consistently outperform co-kriging schemes. In contrast to co-kriging, NNs are able to detect nonlinear correlations between fidelity levels more effectively and are capable of dealing with large data sets. In addition, the hyperparameter selection for the proposed multi-fidelity NNs is automatized based on cross validation and Bayesian optimization, and the tuned NN models consistently yield very low prediction errors. 

The proposed schemes have been tested on a series of manufactured benchmarks and applied to a parametrized PDE problem. In the latter, the goal is to evaluate an output functional of the PDE solution that features an oscillating input-output dependence, and the LF model is constructed through the reduced basis method while the HF model is given by detailed finite element analysis. Numerical results show that the accuracy of the predictions through  MF regression is mainly driven by the reliability of the reduced order model, rather than the amount of HF data fed into the NNs.

A promising direction for future work is to inspect whether changing the value of the balance parameter $\alpha$ during the training process can improve the model accuracy. As the ratio $\mathrm{MSE}_\texttt{HF}/\mathrm{MSE}_\texttt{LF}$ evolves during the training process, a corresponding adaptation of $\alpha$ could be a reasonable improvement of the current schemes.

\bibliographystyle{abbrv}
\bibliography{references.bib}

\setcounter{table}{0}
\renewcommand{\thetable}{A.\arabic{table}}

\section*{Appendix A: hyperparameter summary}

Tables \ref{summary} and \ref{tab:helmholtz} show the values of NN hyperparameters resulting from the HPO in different numerical examples. We recall that $\alpha$ is the parameter in all-in-one networks that balances the fidelity levels' contributions to training error, $\lambda$ is a $L^2$-regularization parameter, and $\eta$ is the learning rate. \textit{Nodes} indicates the number of nodes in each layer whose size is left to be determined through the HPO. Although NN hyperparameters should generally not be analyzed separately, the following comments about the individual hyperparameters can be made:
\begin{itemize}
	\item In all test cases, the optimization procedure chooses $\alpha < 0.1$, which corresponds to weighting the LF data one order of magnitude more than their HF counterpart.
	\item The learning rates $\eta$ in models belonging to the same NNMFR class (all-in-one or multilevel) are generally at the same order of magnitude.
	\item In most cases, the AdaMax algorithm is the best performing optimizer.
	\item In the benchmarks, Glorot uniform weight initialization is preferred, whereas standard uniform initialization performs best when predicting the IRI in the acoustic horn problem.
\end{itemize}

{\small{
	\begin{table}[h!]
		\caption{Hyperparameter values of the NNMFR used for different artificial test cases. Benchmark 1 involves linear correlation, benchmark 2 discontinuous functions, and benchmark 3 non-linear correlation.}
		\renewcommand{\arraystretch}{1.3}
		\centering%
		\begin{tabular}{@{}lccccccc@{}}
			\hline
			Test & Model &Weight initializer &Depth $\times$ width & $\alpha$ & $\lambda$ &Optimizer &$\eta$ \\
			\hline
			\multirow{3}*{1} &  &  & & &\\
			& Intermediate & Glorot normal & 2$\times$59 &4.59$\times 10^{-4}$ & 2.28$\times 10^{-3}$ & Adam & 6.35$\times 10^{-3}$\\
			& GPmimic &Glorot uniform & 3$\times$33 &5.01$\times 10^{-4}$ &1.16$\times 10^{-4}$ &AdaMax &3.90$\times 10^{-3}$\\
			& 2-step &uniform & 1$\times$34 &- &1.01$\times 10^{-4}$ &AdaMax &5.17$\times 10^{-2}$\\
			& 3-step-lin &Glorot uniform &1$\times38$ & -  & 1.08$\times10^{-7}$  &Adam &2.06$\times10^{-2}$ \\
			\multirow{3}*{2} &  &  & & &\\
			& Intermediate & uniform & $2\times30$ &9.22$\times10^{-3}$ &1.94$\times10^{-3}$ &AdaMax &3.41$\times10^{-2}$\\
			& GPmimic & Glorot uniform & $4\times10$ &2.71$\times10^{-2}$ &2.17$\times10^{-4}$ &AdaMax &4.42$\times10^{-2}$\\
			& 2-step & normal & 1$\times60$ & - & 1.01$\times10^{-3}$ & AdaMax & 3.65$\times10^{-3}$\\
			& 3-step & uniform & 1$\times98$ & - & 2.30$\times10^{-4}$ & AdaMax & 1.02$\times10^{-3}$ \\
			\multirow{3}*{3} \\
			& Intermediate &Glorot uniform & $3\times39$ &4.02$\times10^{-4}$ &2.85$\times10^{-5}$
			&Adam &1.55$\times10^{-2}$\\
			& GPmimic & Glorot normal & $4\times22$ &9.74$\times10^{-1}$ &0 &AdaMax &1.14$\times10^{-4}$\\
			& 2-step & uniform & 1$\times44$ & - & 1.02$\times10^{-4}$ & AdaMax & 4.84$\times10^{-3}$\\
			& 3-step & Glorot uniform & 1$\times62$ & - & 2.35$\times10^{-4}$ & Adam & 1.25$\times10^{-4}$\\
			\hline
		\end{tabular}
		\label{summary}
	\end{table}
}}

{\small{
\begin{table}[h!]
 \caption{Hyperparameter values of the considered NN architectures regarding the problems in Section 5; here \textit{P} indicates the number of parameters considered.}
\renewcommand{\arraystretch}{1.3}
   \begin{tabular}{l*{9}{c}c}
    \hline
    P & Model & Weight initializer & $\alpha$ & $\lambda$ & Optimizer & Nodes & $\eta$ & Elapsed time (s)\\
    \hline
    1 & Inter. & Uniform & 3.19$\times 10^{-2}$ & 2.65$\times 10^{-3}$ & Adam & 114 & 4.27$\times 10^{-4}$ & 1426 \\ 
    1  & 2-step & Uniform & - & 1.21$\times 10^{-4}$ & Adamax & 18 & 1.76$\times 10^{-3}$ & 997 \\
    1  & 3-step & Glorot unif. & - & 2.01$\times 10^{-4}$ & Adam & 6 & 1.42$\times 10^{-3}$ & 995\\
    5  & 3-step & Uniform & - & 7.33$\times 10^{-2}$ & Adam & 48 & 7.57$\times 10^{-4}$ & 1089\\
    \hline
    \end{tabular} 
   \label{tab:helmholtz}
\end{table}
}}

\section*{Appendix B: HF and LF models for the acoustic horn problem}

To derive the HF model related to the test case of Section 5, we employ the Galerkin-finite element method. We write the weak formulation of problem \eqref{eq: Helmholtz_system}: given $\boldsymbol{\mu} \in \mathcal{D}$, find $p(\boldsymbol{\mu}) \in V$ s.t.
\begin{equation}
    a(p(\boldsymbol{\mu},u;\boldsymbol{\mu}) = f(u; \boldsymbol{\mu}) \qquad \forall u \in V 
    \label{eq: weak_form}
\end{equation} 
where $V = H^{1}(\Omega(\boldsymbol{\mu}_g)) = \{q \in L^2(\Omega(\boldsymbol{\mu}_g)) \, : \,  \partial{q}/\partial{x_j} \in L^2(\Omega(\boldsymbol{\mu}_g)), \: j\in \{1,2\} \}$ and the bilinear form $a(\cdot,\cdot,\boldsymbol{\mu}):V\times V \rightarrow \mathbb{C}$ and the linear form $f(\cdot;\boldsymbol{\mu}):V \rightarrow \mathbb{C}$ are defined, respectively, by
\begin{align}{2}
    &a(p,u;\boldsymbol{\mu}) = \int_{\Omega(\boldsymbol{\mu}_g)} \{ \nabla p \cdot \nabla \bar{u} - k^2 p \bar{u} \} \,d\Omega + ik\int_{\Gamma_o \cup \Gamma_i} p \bar{u}\,d\Gamma + \frac{1}{2R}\int_{\Gamma_o} p \bar{u}\,d\Gamma \\
    &f(u;\boldsymbol{\mu}) = 2ikA\int_{\Gamma_i} \bar{u}\,d\Gamma
\end{align} 
Next, we introduce a conforming triangulation $\mathcal{T}_h = \{\Delta_k\}_{k=1}^{n_e}$ of the domain $\Omega$ and seek a HF approximation $p_h(\boldsymbol{\mu}) \in V_h$ as a globally continuous, piecewise linear, function belonging to a finite-dimensional space $V_h \subset V$. In our case, $V_h$ is spanned by a set of basis functions $\{\phi_i\}_{i=1}^n$ for the space $V_h$ consisting of a set of $n$ piecewise polynomial nodal basis functions on $\mathcal{T}_h$. The Galerkin-finite element approximation of \eqref{eq: weak_form} thus results in the following $n$ dimensional linear system:
\begin{equation}
   \textbf{A}(\boldsymbol{\mu})\textbf{p}_h (\boldsymbol{\mu}) = \textbf{f}(\boldsymbol{\mu})
   \label{eq: horn_lin_sys2}
\end{equation}
where $\textbf{A}_{ij}(\boldsymbol{\mu}) = a(\phi_j, \phi_i; \boldsymbol{\mu}), \, \textbf{f}_i(\mu) = f(\phi_i;\boldsymbol{\mu}), \, \text{for} \leq i,j \leq n$ and $\textbf{p}_h$ is the vector of coefficients $\{p_i\}_{i=1}^{n}$ such that the projection of $p$ onto $V_h$ is $p_h(\mathbf{x})= \sum_{i=1}^n p_i \phi_i(\mathbf{x})$. This latter formula allows us to state a one-to-one correspondence between the finite element functions $p_h(\boldsymbol{\mu}) \in V_h$ and their discrete counterparts ${\bf p}_h(\boldsymbol{\mu}) \in \mathbb{R}^{N_h}$.

To derive the LF model, we employ the reduced basis (RB) method \cite{HRS,QMN}, which is briefly recalled here. The RB method is a projection-based reduced order modeling technique, addressing the repeated solution of parametrized PDEs, which allows to dramatically reduce the dimension of the discrete problems arising from numerical approximation. The strategy adopted in RB methods consists in the projection of the HF problem upon a subspace made of specially selected basis functions, built from a set of HF solutions corresponding to suitably chosen parameters (or snapshots), e.g., through proper orthogonal decomposition (POD). Later, a (Petrov-)Galerkin projection onto the RB space is employed to generate the ROM. 

Starting from the FOM \eqref{eq: horn_lin_sys2}, i.e. find $\textbf{p}_h(\boldsymbol{\mu})$ such that $\textbf{A}(\mu)\textbf{p}_h(\boldsymbol{\mu}) = \textbf{f}(\mu)$, the idea of a projection-based ROM is to approximate $\textbf{p}_h(\mu) \approx \textbf{V}\textbf{p}_N(\mu)$ as a linear combination of basis functions, for a vector of unknown reduced degrees of freedom $\textbf{p}_N(\mu)$ of reduced dimension $N \ll n$. This latter is sought by imposing that 
\[ \textbf{W}^\text{T}(\textbf{A}(\boldsymbol{\mu})\textbf{V}\textbf{p}_N(\boldsymbol{\mu}) - \textbf{f}(\boldsymbol{\mu})) = {\bf 0},\]
a condition which enforces the orthogonality of the residual to a subspace spanned by a suitable test basis $\textbf{W}\in \mathbb{R}^{n\times N}$. The ROM reads: find $\textbf{p}_N \in \mathbb{R}^N$ such that
\begin{equation}
    {\textbf{A}_N(\boldsymbol{\mu})\textbf{p}_N(\boldsymbol{\mu}) = \textbf{f}_N(\mu)}
    \label{eq: ROM}
\end{equation}
where $\textbf{A}_N(\boldsymbol{\mu} )= \textbf{W}^\text{T} \textbf{A}(\boldsymbol{\mu}) \textbf{V}$ and $\textbf{f}_N(\boldsymbol{\mu}) = \textbf{W}^\text{T} \textbf{f}(\boldsymbol{\mu} )$. A Galerkin projection results if $\textbf{W} = \textbf{V}$. As for the HF approximation, the RB approximation reflects a one-to-one correspondence between the function $p_N^h(\boldsymbol{\mu}) \in V_h$ and its finite-dimensional counterpart  $\textbf{V}{\bf p}_N(\boldsymbol{\mu}) \in \mathbb{R}^{N_h}$.

Although the dimension of the RB problem \eqref{eq: ROM} is very small as compared to the FOM \eqref{eq: horn_lin_sys2}, the assembling of the former system still depends in general on the dimension $n$ of the HF system for any $\boldsymbol{\mu}  \in \mathcal{D}$. A convenient situation arises when the HF arrays in \eqref{eq: horn_lin_sys2} can be written as
\[
\textbf{A}(\boldsymbol{\mu})  = \sum_{q=1}^{Q_A} \Theta_q^a(\boldsymbol{\mu}) \textbf{A}_q, \qquad
\textbf{f}(\boldsymbol{\mu})  = \sum_{q=1}^{Q_f} \Theta_q^f(\boldsymbol{\mu}) \textbf{f}_q.
\]
By virtue of this property – which we refer to as affine parametric dependence of $\textbf{A}(\boldsymbol{\mu})$ and $\textbf{f}(\boldsymbol{\mu})$, the assembling of the system \eqref{eq: ROM} during the online stage can be made efficient, since the arrays $\textbf{W}^\text{T} \textbf{A}_q \textbf{V}$,  $q=1,\ldots,Q_A$ and $\textbf{W}^\text{T} \textbf{f}_q$,  $q=1,\ldots,Q_F$, can be pre-computed and stored during a possibly expensive offline stage.

Since in the Helmholtz problem with $p>1$ parameters we deal with parametrized shape deformations, the FOM arrays $\textbf{A}(\boldsymbol{\mu})$ and $\textbf{f}(\boldsymbol{\mu})$ are nonaffine functions of $\boldsymbol{\mu}$, so we employ hyper-reduction through the discrete empirical interpolation method (DEIM) \cite{DEIM} and its matrix version (MDEIM) \cite{MDEIM,MDEIM1} to compute approximate affine decompositions as \cite{Horn_Manzoni}
\begin{equation*}
    \textbf{f}(\boldsymbol{\mu}) \approx \textbf{f}_m(\boldsymbol{\mu})=\sum_{k=1}^{M_\text{f}}\theta_k^{\text{f}}(\boldsymbol{\mu})\textbf{f}_k\,, \qquad
    \textbf{A}(\boldsymbol{\mu})  \approx \textbf{A}_m(\boldsymbol{\mu}) = \sum_{k=1}^{M_A} \Theta_k^a(\boldsymbol{\mu}) \textbf{A}_k\,,
    \end{equation*}
where $\textbf{f}_k$, $k = 1,\ldots, M_{\text{f}}$, and $\textbf{A}_k$, $k=1,\ldots,M_A$  
are precomputable vectors and matrices, respectively, and independent of $\boldsymbol{\mu}$. In this way, we can approximate the ROM arrays as 
\begin{equation*}
    \textbf{f}_N(\boldsymbol{\mu})\approx \textbf{f}^m_N(\boldsymbol{\mu})=\sum_{k=1}^{M_\text{f}}\theta_k^{\text{f}}(\boldsymbol{\mu})\textbf{f}^k_N\,, \qquad 
        \textbf{A}_N(\boldsymbol{\mu})\approx \textbf{A}^m_N(\boldsymbol{\mu})=\sum_{k=1}^{M_\text{a}}\theta_k^{\text{a}}(\boldsymbol{\mu})\textbf{A}^k_N\,,
\end{equation*}
where $\textbf{f}^k_N=\textbf{W}^\text{T} \textbf{f}_k \in \mathbb{R}^N$, $k = 1,\ldots, M_{\text{f}}$, and $\textbf{A}^k_N=\textbf{W}^\text{T} \textbf{A}^k \textbf{V}  \in \mathbb{R}^{N \times N}$, $k=1,\ldots,M_A$.
Taking advantage of hyper-reduction, we recover the following hyper reduced order model: find $\textbf{p}^m_N(\boldsymbol{\mu}) \in \mathbb{R}^N$ s.t.
\begin{equation}
     \textbf{A}^m_N(\boldsymbol{\mu})\textbf{p}^m_N (\boldsymbol{\mu}) = \textbf{f}^m_N(\boldsymbol{\mu}).
    \label{eq: hyper_form}
\end{equation} 
Due to its small dimension, the solution of the system \eqref{eq: hyper_form} can be very fast and computationally inexpensive, allowing us to generate many instances of the output of interest, which can be used as LF training data. Finally, we denote by $p_N^m(\boldsymbol{\mu}) \in V_h$ the finite element approximation of the problem corresponding to the vector $\mathbf{V} \textbf{p}^m_N (\boldsymbol{\mu}) \in \mathbb{R}^{N_h}$.

\end{document}